\documentclass[12pt,a4paper]{article}  

\usepackage{amsmath}
\usepackage{amssymb}
\usepackage{amsthm}
\usepackage{tikz-cd}
\usepackage{tensor}
\usepackage{mathtools}   
\usepackage{hyperref} 
\usepackage{bbm}

\theoremstyle{plain}
\newtheorem{theorem}{Theorem}
\newtheorem{lemma}[theorem]{Lemma}
\newtheorem{proposition}[theorem]{Proposition}
\newtheorem{corollary}[theorem]{Corollary}
\theoremstyle{definition}
\newtheorem{definition}[theorem]{Definition}
\theoremstyle{remark}
\newtheorem{remark}[theorem]{Remark}

\def\AA            {B}
\def\actle         {\blacktriangleright}
\def\actlE         {\,{\actle}\,}
\def\actre         {\vartriangleright}

\def\actrE         {\,{\actre}\,}
\def\alph          {\alpha}  
\def\be            {\begin{equation}}
\def\bimod         {\text-{\rm bimod}}
\def\Bmod          {\,{\mathrm{mod}}\,}

\def\calc          {{\mathcal C}}
\def\calm          {{\mathcal M}}
\def\caln          {{\mathcal N}}
\def\cir           {\,{\circ}\,}
\def\coev          {\mathrm{coev}}
\def\coHom         {\mathrm{coHom}}
\def\Colon         {:\quad}

\def\Cong          {\,{\cong}\,}
\newcommand\deltal[3]{{}_{}^{\ihl}\delta^{#2}_{#1,#3}}
\def\deltaL        {{}_{}^{\ihl}\delta}
\newcommand\deltar[3]{{}_{}^{\ihr}\delta^{#2}_{#1,#3}}
\def\deltaR        {{}_{}^{\ihr}\delta}
\def\ee            {\end{equation}}
\def\End           {\mathrm{End}}
\def\Enumerate     {\def\leftmargini{1.34em}~\\[-1.42em]\begin{enumerate}}

\def\eq            {\,{=}\,}
\newcommand\eqc[1] {\!\!\stackrel{\eqref{#1}}\cong\!\!}
\def\ev            {\mathrm{ev}}

\def\Hom           {\mathrm{Hom}}

\def\icH           {\underline{\coHom}}
\def\icoHom        {\icH}
\def\icoHomr       {\icH^{\ihr}}
\def\icoHoml       {\icH^{\ihl}}

\def\id            {\mathrm{id}}
\def\Id            {\mathrm{id}}

\def\iev           {\underline{\mathrm{ev}}}

\def\iH            {\underline{\Hom}}
\def\ihl           {{\mathrm l}}

\def\iHomr         {\iH^{\ihr}}
\def\iHoml         {\iH^{\ihl}}
\def\ihr           {{\mathrm r}}
\def\imu           {\underline{\mu}}
\def\iN            {\,{\in}\,}
\def\inv           {^{-1}}
\def\ko            {{\ensuremath{\Bbbk}}}

\def\moD           {{\rm mod}\text-}
\def\Mod           {\text-{\rm mod}}
\newcommand\nxl[1] {\\[#1mm]}
\newcommand\Nxl[1] {\\[-1.3em]\\[#1mm]}
\def\oat           {\otimes^{A}}
\def\oaT           {\,{\oat}\,}
\def\one           {\TI}
\def\opp           {^{\rm opp}}
\def\ota           {\otimes_{\!A}}
\def\otA           {\,{\ota}\,}
\def\otAz          {\,{\otimes_{\!A_2}}\,}
\def\oti           {\,{\otimes}\,}
\def\otik          {\,{\otimes_\ko}\,}

\newcommand\rarr[1]{\xrightarrow{\,{#1}\,}}
\newcommand\Rarr[1]{\,{\xrightarrow{\,{#1}\,}}\,}
\def\RL            {Q}  
\def\TI            {1}
\def\Times         {\,{\times}\,}
\def\To            {\,{\to}\,}

\def\Tocong        {\Rarr\cong}
\def\tpl           {\bullet}
\def\tpL           {\,{\tpl}\,}
\def\tpr           {\otimes}
\def\tpR           {\,{\tpr}\,}
\def\vect          {\mathrm{vect}}

\makeatletter
\newcommand*{\relrelbarsep}{.386ex}
\newcommand*{\relrelbar}{%
  \mathrel{%
      \mathpalette\@relrelbar\relrelbarsep }}
\newcommand*{\@relrelbar}[2]{%
  \raise#2\hbox to 0pt{$\m@th#1\relbar$\hss}%
    \lower#2\hbox{$\m@th#1\relbar$}}
\providecommand*{\rightrightarrowsfill@}{%
      \arrowfill@\relrelbar\relrelbar\rightrightarrows }
\providecommand*{\leftleftarrowsfill@}{%
       \arrowfill@\leftleftarrows\relrelbar\relrelbar }
\providecommand*{\xrightrightarrows}[2][]{%
          \ext@arrow 0359\rightrightarrowsfill@{#1}{#2}}
\providecommand*{\xleftleftarrows}[2][]{%
    \ext@arrow 3095\leftleftarrowsfill@{#1}{#2}}
\makeatother

\begin{document}

\numberwithin{equation}{section}
\numberwithin{theorem}{section}


  \thispagestyle{empty}
  \begin{flushright}
     {\sf ZMP-HH/23-11}\\
     {\sf Hamburger$\;$Beitr\"age$\;$zur$\;$Mathematik$\;$Nr.$\;$945}
     \\[2mm] June 2023
  \end{flushright}

	     
\vskip 2.7em

\begin{center}     
        {\bf \Large Grothendieck-Verdier duality in categories of bimodules
        \\[4pt] and weak module functors}

\vskip 2.6em

{\large 
J\"urgen Fuchs\,$^{a},\,$
Gregor Schaumann\,$^{b},\,$
Christoph Schweigert\,$^{c},\,$
Simon Wood\,$^{d,c}$
}

\vskip 15mm

 \it$^a$           
 Teoretisk fysik, \ Karlstads Universitet \\
 Universitetsgatan 21, \  S\,--\,651\,88\, Karlstad
 \\[9pt]
 \it$^b$
 Mathematische Physik, \ Institut f\"ur Mathematik, \ Universit\"at W\"urzburg \\
 Emil-Fi\-scher-Stra\ss e 31, \ D\,--\,97\,074\, W\"urzburg
 \\[9pt]
 \it$^c$
 Fachbereich Mathematik, \ Universit\"at Hamburg\\
 Bereich Algebra und Zahlentheorie\\
 Bundesstra\ss e 55, \ D\,--\,20\,146\, Hamburg
 \\[9pt]
 \it$^d$
 School of Mathematics, \ Cardiff University\\
 Abacws, \ Senghennydd Road, \ Cardiff, CF24 4AG

\end{center}

\vskip 4.7em

\noindent{\sc Abstract}\\[3pt]
Various monoidal categories, including suitable representation categories
of vertex operator algebras, admit natural Grothen\-dieck-Verdier duality 
structures. We recall that such a Grothen\-dieck-Verdier category comes with 
two tensor products which should be related by distributors obeying pentagon 
identities. We discuss in which circumstances these distributors are isomorphisms. 
This is achieved by taking the perspective of module categories over monoidal 
categories, using in particular the natural weak module functor structure
of internal Homs and internal coHoms. As an illustration, we exhibit these
concepts concretely in the case of categories of bimodules over associative
algebras.

\newpage

\section{Introduction}

Dualities play a pivotal role in various applications of monoidal 
categories. A familiar example of a duality is a rigid structure. For many
applications, for instance in representation theory and in linear logic, it is,
however, necessary to consider more general notions of duality. For example, if 
a tensor product is not exact, then the category does not admit a rigid duality. 

A framework that generalizes rigidity and in which structures familiar from 
rigidity, such as internal Homs, persist, is furnished by a $\ast$-autonomous 
structure, also known as a \emph{Grothen\-dieck-Ver\-dier duality}. Categories 
endowed with a Grothendieck-Verdier duality have a wide range of applications, 
from linear logic to the representation theory of vertex operator algebras.
In this paper we present some important aspects of Grothendieck-Verdier dualities,
concentrating on such which are related to internal Hom and coHom functors.
Our point of view is the one of module categories. In particular, we will regard
a monoidal category as a module category over itself. Specifically, a 
crucial ingredient is the structure 
of internal Homs and coHoms as weak module functors. This allows us to introduce 
distributors, which play a role similar to the associator for the tensor product.
Before discussing these concepts in appropriate generality, we illustrate some
of them in the, arguably, simplest relevant situation -- finite-dimensional
bimodules over finite-dimensional algebras. After presenting the general theory,
we return to this example and present explicit formulas for the distributors
in categories of bimodules.


\section{Bimodules} \label{sec:bimod}

To start with, let us consider some very basic mathematical objects which are
familiar already to undergraduates. We fix an algebraically closed field \ko. 

Let $A$ be a (unital associative)
\ko-algebra. The category of $A$-bimodules can be endowed with a monoidal structure 
given by the tensor product $\ota$ over $A$, which has $A$ as 
its monoidal unit. The tensor product $\ota$ is defined by the coequalizer diagram
  \be
  B_1\otik A\otik B_2 \rarr{\,\varphi\,} B_1\otik B_2\rarr{~\,~} B_1 \otA B_2 \,.
  \ee
Here the first arrow is the linear map
  \be
  \varphi = \rho^{B_1}\otik \id_{B_2} - \id_{B_1}\otik \lambda^{B_2} ,
  \ee
where $\rho^{B_1}$ is the right $A$-action on $B_1$ and $\lambda^{B_2}$ the left
$A$-action on $B_2$. Here, as well as largely below, we suppress the associator
$\alph$ of a monoidal category, and likewise we will suppress the unitors $l$ and $r$.

The tensor product $\ota$ is right exact.
Since our primary interest in this note is in duality structures, it is convenient
to impose suitable finiteness conditions. Specifically, unless stated otherwise, in 
the sequel all $\ko$-algebras are assumed to be finite-dimensional, and all modules
and bimodules to be finite-dimensional as \ko-vector spaces. Then in particular the
tensor product is right exact and equips the category $A\bimod$ of $A$-bimodules with
a monoidal structure, and we can make
contact to the theory of coalgebras and comodules.

When $A$ is commutative, then every left $A$-module is canonically a right module 
and even a bimodule. Accordingly we may then also consider the category $A\Mod$ of
left $A$-modules, which is a monoidal subcategory of the category $A\bimod$ of
$A$-bimodules, with $A$ as a left module as the monoidal unit. One motivation to
consider commutative algebras in this elementary study comes from the representation
theory of vertex operator algebras. Indeed, commutative algebras can be seen as 
particular examples of vertex operator algebras \cite[Sect.\,1.3]{FRbe2}, and 
finite-dimensional
algebras as particular examples of $C_2$-cofinite vertex operator algebras.
Vertex operator algebras that are commutative algebras are known to arise as 
the result of BRST cohomology or quantum Hamiltonian reduction, see e.g.\ Section
2.2 of \cite{lizu4} for an infinite-dimensional example. BRST cohomology can also
result in finite-di\-mensional commutative algebras.\,\footnote{~
 We thank Sven M\"oller for pointing out an (unpublished) example of such a
 commutative algebra.}

A fundamental operation in linear algebra is to take the linear dual of a vector space.
The category $\vect_\ko$ of finite-dimensional vector spaces is \emph{rigid}, that
is, every vector space $V$ has both a left dual $^\vee V$ and a right dual $V^\vee$ 
(both of which are isomorphic to $\Hom_\ko(V,\ko)$ since $\vect_\ko$ is symmetric)
endowed with evaluation and coevaluation maps, given by
  \be
  \begin{array}{rlrl}
  \ev \colon~ V^\vee \otik V & \!\!\! \rarr{~\,} \ko \,, \qquad
  & \coev \colon~ \ko & \!\!\! \rarr{~\,} V \otik V^\vee ,
  \\[4pt]
  f \otik v & \!\!\! \longmapsto f(v) \,,
  & \xi & \!\!\! \longmapsto \xi \sum_i v_i \otik v^i
  \end{array}
  \ee
for the right dual, and analogously for the left dual, where $(v_i)$ is a basis
of $V$ and $(v^i)$ the dual basis of $V^\vee_{}\!$.
The notion of rigidity pervades the theory of Hopf algebras and quantum topology.
Indeed, a (quasi-)Hopf algebra is essentially an algebra $H$ endowed with additional
structure that is precisely such that the category $H\Mod$ of finite-dimensional
$H$-modules is rigid monoidal. 

Taking the double dual gives a monoidal endofunctor $(-)^{\vee\vee}$ of the monoidal
category 
$\vect_\ko$. The double dual $V^{\vee\vee}$ of a finite-dimensional vector space is 
canonically identified with $V$. More precisely, the category $\vect_\ko$ is pivotal,
i.e.\ there is a monoidal isomorphism $\Id_{\vect_\ko} {\rarr{\,}}\, (-)^{\vee\vee}$
between the identity functor and the double dual. For an arbitrary rigid monoidal 
category $\calc$ such a monoidal isomorphism may or may not exist; and if it exists,
then there might be several of them. The choice of a specific one is then called 
a pivotal structure on $\calc$. For $\calc \eq H\Mod$ the
category of finite-dimensional modules over a Hopf algebra $H$, a pivotal structure
amounts to a finding a group-like element $g\iN H$ such that the square of the 
antipode $S$ of $H$ is given by $S^2(h) \eq ghg^{-1}$.

 \smallskip

It is tempting to regard rigidity as the generic paradigm for a duality structure
on monoidal categories. However, already long ago more general duality structures
have been encountered (compare e.g.\ Remark \ref{rem:BOur3takeu3} below). Such
generalized dualities occur in important situations, in particular for representation 
categories of vertex operator algebras to which the HLZ-theory of tensor 
products \cite{HLZ} applies. The purpose of the present section is to illustrate
these structures in the simplest possible case: for bimodules over a $\ko$-algebra.

Indeed it is easy to see that in general the category $A\bimod$ cannot be rigid,
even if the algebra $A$ is finite-dimensional and only finite-dimensional
bimodules are considered. In an
abelian rigid monoidal category, the functors $X\oti{-}$ and ${-}\oti X$ of tensoring 
from the left or from the right with an object $X$ are exact \cite[Prop.\,4.2.1]{EGno}.
This property does not hold for bimodules, in general. A simple instructive 
counter example is the following: For any field \ko, the algebra of dual numbers, 
i.e.\ the two-dimensional quotient
  \be
  A_2 := \ko[x]\,/\,\langle x^2\rangle
  \label{eq:def:A2}
  \ee
of the polynomial algebra $\ko[x]$, is a commutative associative unital algebra. 
A canonical basis of $A_2$ is given by the classes $[1] \eq 1\Bmod x^2 $ and 
$[x] \eq x\Bmod x^2$; to unburden notation, we drop the class symbols and just 
write $1$ and $x$ for the elements of this basis. The algebra $A_2$ has, up to
isomorphism, a single simple module $S$, which is one-dimensional with generator 
$s$; the element $x$ of $A_2$ acts on it by zero. The only other indecomposable 
finite-dimensional $A_2$-module, up to isomorphism, is the free module $P \eq A_2$
of rank 1, with basis elements $1$ and $x$; the element $x\iN A_2$ acts as 
$1\,{\mapsto}\, x$ and $x\,{\mapsto}\, x^2 \eq 0$. 

The following argument shows that the simple module $S$ is not flat. The injective 
morphism $\iota\colon S\Rarr~ P$ is given by $s\,{\mapsto}\, x$. One immediately 
verifies that $S\otAz S \,{\cong}\, S$ and, since $P$ is the monoidal unit,
$S\otAz P \,{\cong}\, S$. The structure map of the coequalizer
$S\otik P\Rarr~ S\otAz P$ maps $s\otik x \,{\mapsto}\, 0$. It follows that
tensoring with $S$ maps the monomorphism $\iota$ to the zero morphism: thus indeed
the module $S$ is not flat. As a consequence the category $A_2\Mod$ is not rigid.

Still, even in the non-rigid case we do have dual vector spaces at our disposal. We
account for their presence by introducing, for any pair $A$ and $A'$ of
algebras, the contravariant functor 
  \be
  G_{\!A,A'} \Colon A\mbox-A'\bimod \rarr~ A'\mbox-A\bimod
  \label{eq:def:GAB}
  \ee
that sends a bimodule $B \iN A\mbox-A'\bimod$ to its linear dual 
$B^\ast \eq \Hom_\ko(B,\ko)$ endowed with the canonical action 
  \be
  (a'.\beta.a)(b) := \beta(a.b.a')
  \label{eq:dualbimod}
  \ee
for $a\iN A$, $a'\iN A'$, $b\iN B$ and $\beta \iN B^\ast$. Obviously this functor
squares to the identity, $G_{A',A} \cir G_{A,A'} 
           $\linebreak[0]$
           {=}\, \Id$. In particular, for $A' \eq A$ we get this way an $A$-bimodule
$A' \eq A$ we get this way an $A$-bimodule $A^\ast$; this furnishes a distinguished 
object in the monoidal category $A\bimod$.

The distinguished bimodule $A^\ast$ is, in general, not isomorphic to the monoidal 
unit $A$ of $A\bimod$, as is already illustrated by the simple example 
  \be
  \AA_3 := \ko[x,y] \,/\, {\langle x^2,y^2,xy\rangle}
  \label{eq:def:A3}
  \ee
of a three-dimensional commutative algebra. A canonical basis of $\AA_3$ is given by 
the classes $[1]$, $[x]$ and $[y]$, which again we abbreviate just as $1$, $x$ and $y$.
The
element $1$ is the unit of $\AA_3$, while $x$ and $y$ satisfy $x x \eq y y \eq x y \eq 0$.
The algebra $\AA_3$ has a unique (up to isomorphism) simple (bi)module $S$, which is
one-dimensional. There is a short non-split exact sequence
  \be
  0 \rarr~ S \,{\oplus}\, S \rarr~ \AA_3 \rarr~ S \,,
  \ee   
for $\AA_3$, from which by dualizing we obtain a short exact sequence
  \be
  0 \rarr~ S \rarr~ \AA_3^\ast \rarr~ S \,{\oplus}\, S \rarr~ 0
  \ee	
for $\AA_3^\ast$. In particular, the modules $\AA_3$ and $\AA_3^\ast$ have
different socle and hence are not isomorphic. (As a consequence, $\AA_3$ does
not admit the structure of a Frobenius algebra; it is in fact the simplest
example of an algebra with this property.)

These observations raise the question: What is the appropriate \emph{categorical}
duality structure (which cannot be rigidity) on the monoidal category $A\bimod$
that accomodates the distinguished bimodule $\AA_3^\ast$ and the duality
\eqref{eq:def:GAB} inherited from the duality on vector spaces?

 \medskip

Before discussing the appropriate duality structure in detail, let us recall 
a few further properties that the monoidal category $A\bimod$ of bimodules over
a finite-dimensional algebra inherits from $\vect_\ko$. These will eventually
find their explanation in the framework of the general duality structures.

We first note the somewhat less well-known fact that for a finite-dimensional 
algebra $A$ the abelian category $A\bimod$ admits a second monoidal product. 
To see this, note that the vector space $A^\ast$ has a natural structure of a 
coalgebra. Further, any right $A$-module $(M,\rho\colon M\otik A\,{\to}\, M)$
becomes a right $A^\ast$-comodule as follows: Select a basis $(a_i)$ of $A$,
which yields a dual basis $(a^i)$ of $A^*$. Then the prescription
  \be
  \tilde\rho(m):= \sum_i m.a_i \otik a^i 
  \ee
defines a right $A^\ast_{}$-coaction. A similar construction turns a left $A$-module
into a left $A^\ast_{}$-comodule. One then defines a tensor product $\otimes^A$ of
comodules as an equalizer
  \be
  B_1\otimes^A B_2 \rarr~ B_1\otik B_2 \rarr\varphi B_1 \otik A^\ast \otik B_2 
  \ee
with $\varphi$ the linear map
  \be
  \varphi = \tilde\rho^{B_1}\otik \id_{B_2} - \id_{B_1}\otik \tilde \lambda^{B_2} .
  \ee
Since this tensor product is defined as a limit, it is, as usual for coalgebras,
left exact. Its monoidal unit is the bimodule $A^\ast\!$.

We thus have two tensor products on $A\bimod$; the tensor product $\otimes^A$ is 
left exact, while $\ota$ is right exact. This suggests to study their left and
right adjoints, respectively. Doing so leads to the notions of internal Homs and 
internal coHoms; we introduce these for general linear monoidal categories:

\begin{definition} \label{def:iHom}
Let $\calc$ be a \ko-linear monoidal category with tensor product $\otimes$,
and let $Y,Z\iN\calc$. The
\emph{internal Hom} $\iHomr(Y,Z)$ is an object representing the functor
$\calc\opp\Rarr~\vect_\ko$ that is given by $X\,{\mapsto}\, \Hom(X\tpR Y,Z)$, 
i.e.\ we have functorial isomorphisms
  \be
  \Hom(X\tpR Y,Z) \cong \Hom(X, \iHomr(Y,Z)) \,.
  \label{eq:def:iHomr}
  \ee
Similarly a second internal Hom $\iHoml$ is obtained by keeping the second tensor 
factor \cite[Remarks\,2.1(3)]{boDr}, i.e.\
  \be
  \Hom(X\tpR Y,Z) \cong \Hom(Y, \iHoml(X,Z)) \,.
  \label{eq:def:iHoml}
  \ee
The category $\calc$ is called \emph{right closed} if all internal Homs $\iHomr(Y,Z)$
exist, and it is called \emph{left closed} if all internal Homs $\iHoml(Y,Z)$ exist.
\end{definition}

A necessary condition for the existence of the internal Homs is 
that the tensor product $\otimes$ is right exact.
The analogous notion for the left exact tensor product $\otimes^A$ on $A\bimod$
is the left adjoint of $\otimes^A$. Again this can be generalized:

\begin{definition} \label{def:icoHom}
Let $\calc$ be a \ko-linear monoidal category with tensor product $\otimes$,
and let $Y,Z\iN\calc$. Then the
\emph{internal coHom} is the object characterized by functorial isomorphisms
  \be
  \Hom(X,Y\oti Z) \cong \Hom(\icoHomr (Z,X),Y) \,.
  \label{eq:def:icoHom}
  \ee
Similarly there is a second internal coHom $\icoHoml$.
\end{definition}

For the case of the monoidal category of $A$-bimodules considered here, internal 
Homs and coHoms do exist. Using the canonical pivotal structure of 
finite-di\-men\-si\-onal vector spaces, they can be expressed
in terms of the tensor products $\otimes^A$ and $\ota$, respectively:
  \be
  ~\iHomr(X,Y) = Y \,{\otimes^A}\, X^\ast \quad \text{and}\quad
  \icoHomr(X,Y) = Y \otA X^\ast ,
  \label{eq:iHom/icoHom.as.tpr/tpl}
  \ee
as well as
  \be
  ~\iHoml(X,Y) = X^\ast \,{\otimes^A}\, Y \quad \text{and}\quad
  \icoHoml(X,Y) = X^\ast \otA Y .
  \label{eq:iHoml/icoHoml.as.tpr/tpl}
  \ee

To show that the category $A\bimod$ does admit internal Homs, we provide them 
explicitly. In fact, for any triple of algebras $A$, $B$ and $C$ and bimodules
${}_{C}M_{B}$, ${}_{B}N_{A}$ and ${}_{C}X_{A}$ there is the adjunction
  \be
  \begin{aligned}
  \Hom_{C,A}({}_{C}M_{B} \, {\otimes_B}\, {}_{B}N_{A},{}_{C}X_{A})
  & \cong \Hom_{C,B}({}_{C}M_{B}, \Hom_{A}(N_{A},X_{A}))
  \nxl1
  & \cong \Hom_{C,B}({}_{C}M_{B}, X_A \oaT (N_A)^{*}_{}) \,,
  \end{aligned}
  \label{eq:Adj-bimod}
  \ee
where we indicate only one of the algebra actions on a bimodule in cases when
it is more relevant than the other, e.g.\ $N_{A}$ is the bimodule ${}_{B}N_{A}$ seen
as a right $A$-module, and where 
in the second expression,
$\Hom_{A}(N_{A},X_{A})$ is regarded as $(C,B)$-bimodule
via the left module structures of $N_{A}$ and $X_{A}$. We can thus read off that
  \be
  \iHomr(N,X) = \Hom_{A}(N_{A},X_{\!A}) \,.
  \label{eq:iHomrNX}
  \ee
Similarly we find
  \be
  \iHoml(N,X) = \Hom_{A}({}_{A}N,{}_{A\!}X) \,.
  \label{eq:iHomlNX}
  \ee
Note that \eqref{eq:iHomrNX} exhibits the internal Hom as a sub-bimodule
  \be
  \Hom_{A}(N_{A},X_{A}) \subset \Hom_\ko(N,X) 
  \ee
of the bimodule obtained by equating the right $A$-actions on $X$ and $N$,
and analogously for \eqref{eq:iHomlNX}. In contrast, internal coHoms are defined 
as quotients. For finite-dimensional bimodules this fits with the description 
\eqref{eq:iHom/icoHom.as.tpr/tpl} and \eqref{eq:iHoml/icoHoml.as.tpr/tpl} of 
internal Homs and coHoms in terms of the two tensor products $\otimes^A$ and 
$\ota$. Also note that according to \eqref{eq:iHomrNX} and
\eqref{eq:iHomlNX}, left and right internal Homs are, in general, non-isomorphic.

Finally we specialize \eqref{eq:Adj-bimod} to the case that $X \eq A^\ast$. 
Restricting further to the special case that $M$ and $N$ are $A$-bimodules we then get
  \be
  \begin{aligned}
  \Hom_{A,A}(M \otA N, G({}_{A}A_{A}))
  & \eqc{eq:Adj-bimod} \Hom_{A,A}(M, \Hom_{A}(N_{A}, (A^\ast)_{A})) 
  \Nxl1
  & ~=~ \Hom_{A,A}(M, \Hom_{\ko,A}(N, \Hom_{\ko}(A,\ko))) 
  \Nxl1
  & \eqc{eq:Adj-bimod} \Hom_{A,A}(M,\Hom_{\Bbbk,\Bbbk}(N \otA A, \Bbbk)) 
  \Nxl1
  & ~\cong~ \Hom_{A,A}(M, G(N)) \,.
  \end{aligned}
  \label{eq:General-G}
  \ee

\begin{remark} \label{rem:infdim}
Recall that to any $A$-bimodule $M$ one can associate its center
  \be
   Z(M) := \{ b\iN M \,|\, a.m \eq m.a ~\text{for\,all}~ a\iN A\}   
  \ee
as a subobject and its trace
  \be
   \mathrm{tr}\, M := M / \langle a.m \,{-}\, m.a \rangle 
  \ee
as a quotient. By comparison with the $A$-bimodule structure on the
vector space $\Hom_\ko(N,X)$ defined above, we see that internal Homs are centers
and internal coHoms are traces.
Also note that for the definition of the center and trace we do not need to require 
that the bimodule $M$ is finite-dimensional. By considering bimodules of the form 
$M \otik N$, this implies that the two tensor products $\otA$ and $\otimes^A$
still exist for infinite-dimensional bimodules.
Similarly, the contravariant functor \eqref{eq:def:GAB} that sends a bimodule to 
its linear dual and the 
adjunction \ref{eq:Adj-bimod} are still present in the infinite-dimensional case.
\end{remark}


\section{Grothendieck-Verdier duality} \label{sec:GV}

\subsection{First definitions} \label{sec:1stdef}

The observations in the preceding section motivate the following

\begin{definition}
Let $\calc \eq (\calc,\tpr,1,\alph,l,r)$ be a monoidal category. 
  \begin{enumerate} \addtolength{\itemsep}{1.8pt}
  \item 
A \emph{dualizing object} 
of $\calc$ is an object $K \iN \calc$ such that for every $y\iN \calc$ the functor 
$x \,{\mapsto}\, \Hom(x
        $\linebreak[0]$
        {\tpr}\, y, K)$ is representable by some object $Gy \iN \calc$
and the so defined contravariant functor $G\colon \calc \To \calc$ is an 
anti-equivalence. 
 \\[2pt]
We thus have isomorphisms
  \be
  \varpi_{x,y}^{} \Colon \Hom(x\tpR y,K) \rarr\cong \Hom(x,Gy) 
  \label{eq:GV1}
  \ee
for $x,y \iN \calc$. $G$ is called the \emph{duality functor with respect to $K$}. 
  \item
A \emph{Grothendieck-Verdier structure} on $\calc$ is the choice of a 
dualizing object $K\iN\calc$.
A \emph{Grothendieck-Verdier category}, or \emph{GV-category}, for short, is a 
monoidal category together with a choice of a Grothendieck-Verdier structure.
\end{enumerate}
\end{definition}

While the functor $G$ depends on the choice of $K$, we suppress it in the notation.\,%
 \footnote{~It is worth noting that this (standard) definition of a GV-category 
 differs from the one in Definition 2.3 of \cite{haLem}, a paper which also 
 discusses the relation to linearly distributive categories.}
This duality structure has been introduced in \cite{BArr} under the name
$\ast$-autono\-mous category, compare also \cite{barr5,barr8,barr9}. The term
Grothendieck-Verdier category dates back to \cite{boDr}. GV-categories appear
for instance in the study of quadratic algebras and operads \cite{mani27} and
also play an important role in linear logic, see e.g.\ \cite{seely3}.

The result \eqref{eq:General-G} in Section \ref{sec:bimod} tells
us that $A\bimod$ is a GV-category with dualizing object $A^*$,
thereby giving a categorical meaning to this distinguished bimodule. This should
not come as a surprise: it is known \cite[Rem.\,3.17]{fScS2} that this fact can
be formulated in a Morita invariant manner. The dualizing object is in fact
\emph{structure}; if $K$ is a dualizing object, then
any object of the form $K\oti D$ with 
$D$ an invertible object of $\calc$ is dualizing as well \cite[Prop.\,1.3]{boDr}.
It should, however, be appreciated that the dualizing object itself need not
be invertible. GV-categories with non-invertible dualizing object indeed arise as
representation categories of vertex operator algebras \cite{garW3}. Below, we
also give an example involving bimodules. 

\begin{definition}
An $r$-\emph{category} is a 
monoidal category for which the monoidal unit $\TI$ is a dualizing object.
\end{definition}

Clearly, every rigid monoidal category is an $r$-category. But there also
exist non-rigid $r$-categories; for an example, see e.g.\ \cite[Ex.\,1.9]{boDr}.

In the case of categories $A\bimod$ of bimodules, certain invertible bimodules can 
be obtained by twisting either the left or right action of $A$ by an algebra 
automorphism $\psi$ of $A$. For the isomorphism class of the bimodule $A_\psi$ 
obtained in this way, $\psi$ only matters up to inner automorphisms. For instance, 
the group of outer automorphisms of the three-dimensional algebra
$\AA_3 \eq \ko[x,y] / {\langle x^2,y^2,xy\rangle}$ considered in \eqref{eq:def:A3} 
is isomorphic to $\mathrm{GL}(2)$, acting on the two generators in the obvious
way. This gives a simple example of a GV-category admitting different GV-structures.

We also have already noticed that in this case the vector space dual $\AA_3^\ast$
is not isomorphic to $\AA_3$ as a bimodule, so $\AA_3$ is not a Frobenius 
algebra. The following consideration shows that $\AA_3^\ast$ is not even
invertible: $\AA_3$ has, up to isomorphism, a single simple module $S$, which is
one-dimensional and on which the elements $x$ and $y$ of $\AA_3$ act as zero. A direct
calculation shows that $\AA_3^\ast\otA \AA_3^\ast \,{\cong}\, S^{\oplus 4}$. Thus if
$\AA_3^\ast$ were invertible, then $S^{\oplus 4}$ would be too, i.e.\ there would 
be a bimodule $X$ such that $S^{\oplus 4}\otA X \,{\cong}\, A$. But this is 
impossible, simply because $A$ is indecomposable as a module over itself.
In fact, one can show that the regular bimodule $\AA_3$ is, up to isomorphism,
the only invertible $\AA_3$-bimodule.


\subsection{Internal Homs}

The existence of a GV-structure on a monoidal category implies much structure
familiar from category theory. In particular, recall from Definition \ref{def:iHom} 
that, for $\calc$ a linear monoidal category and $y,z\iN\calc$, the internal Hom
$\iHomr(y,z) \iN \calc$ is an object representing the functor $x \,{\mapsto}\, 
\Hom(x\tpR y,z)$, as expressed by the functorial isomorphisms \eqref{eq:def:iHomr}.
Since $\calc$ is not assumed to be braided, there is also a separate variant 
$\iHoml$ of the internal Hom in which the 
second tensor factor is kept, with functorial isomorphisms \eqref{eq:def:iHoml}.

A GV-category admits internal Homs. In fact, 
from the defining properties of the functor $G$ it follows that
  \be
  \iHomr(x,z)\cong G(x\tpR G^{-1} z) \qquad\text{and}\qquad
  \iHoml(x,z)\cong G^{-1}(Gz\tpR x) \,.
  \label{eq:iHom=}
  \ee
In particular we have
  \be
  \iHomr(x,K) \cong G(x\tpR G^{-1}(K)) \cong G(x\tpR 1) \cong G(x) \,,
  \ee
i.e.\ there is an isomorphism 
  \be
  \iHomr(-,K) \cong G
  \label{eq:iHom(,K)=G}
  \ee
of functors. This explicitly shows how the duality functor $G$ is determined 
in terms of the dualizing object $K$.

There are canonical evaluation morphisms
  \be
  \iev^{\ihr}_{x,y} \colon~~ \iHomr(x,y) \tpr x \rarr{} y \quad\text{ and } \quad 
  \iev^{\ihl}_{x,y} \colon~~ x \tpR \iHoml(x,y) \rarr{} y
  \label{eq:iev'}
  \ee
for the internal Homs, given by the images of the identity morphism in the spaces
$\End(\iHomr(x,y))$ and $\End(\iHoml(x,y))$ under the defining adjunctions 
\eqref{eq:def:iHomr} and
\eqref{eq:def:iHoml}, respectively. When combined with \eqref{eq:iHom(,K)=G},
this immediately gives left and right evaluation morphisms for the tensor product.
In contrast, compatible coevaluations do not exist, in general; otherwise the
tensor product on any abelian GV-category would be exact.
Given the evaluations, associative multiplication morphisms
  \be
  \begin{aligned}
  \imu^{\ihr}_{x,y,z} \Colon& \iHomr(y,z) \tpR \iHomr(x,y) \rarr{~} \iHomr(x,z)
  \\[2pt]
  \text{and} \qquad
  \imu^{\ihl}_{x,y,z} \Colon& \iHoml(x,y) \tpR \iHoml(y,z) \rarr{~} \iHoml(x,z)
  \end{aligned}
  \label{eq:def:imu}
  \ee
are obtained by standard arguments \cite[Sect.\,7.9]{EGno}. In particular,
$\iHoml(z,z)$ and $\iHomr(z,z)$ have the structure of (unital) associative algebras.


\subsection{A second tensor product} \label{sec:2ndtp}

It is natural to introduce on a GV-category a second monoidal structure by
  \be
  x\tpl y:= G^{-1} (Gy\tpR Gx) \,.
  \label{eq:def:seccirc}
  \ee
The two tensor products are, in general, different, and hence in general
$G(x\tpR y)$ and $G(y)\tpR G(x)$ are not isomorphic, in contrast to
the case of rigid categories. The double dual $G^2$, however, is monoidal
\cite[Prop.\,5.2]{boDr}. It describes how cyclic invariance of invariant tensors 
of Grothendieck-Verdier categories is violated \cite[Rem.\,2.1(2)]{boDr}:
  \be
  \begin{aligned}
  \Hom(G^2(y)\tpR x,K ) & \cong \Hom(G^2(y),Gx)
  \\ 
  & \cong \Hom(x,Gy) \cong \Hom(x\tpR y,K) \,.
  \end{aligned}
  \ee

There is also another variant of the second monoidal structure, given by
  \be
  x\tpl' y := G (G^{-1}y\tpR G^{-1}x) \,.
  \ee
The monoidal structure on the functor $G^2$ provides a canonical identification of
$\tpl$ and $\tpl'$; accordingly we will henceforth identify $\tpl$ and $\tpl'$.

\begin{remark} 
If $\calc$ is an $r$-category, then the monoidal
products $\tpr$ and $\tpl$ are connected by canonical functorial morphisms
  \be
  \varphi_{x,y}^{} \Colon x \tpR y \rarr~ x \tpL y \,.
  \label{eq:X*Y->X.Y}
  \ee
The morphisms $\varphi_{x,y}^{}$ are functorial in $x$ and $y$ and compatible with
the associativity constraints for the two
tensor products. In general they are not isomorphisms.
\end{remark} 

\begin{remark} 
As follows from \eqref{eq:iHom/icoHom.as.tpr/tpl},
in the categories of bimodules considered in Section \ref{sec:bimod}, the two tensor
products $\tpr$ and $\tpl$ are realized as the
\(\ota\)- and \(\otimes^A\)-tensor products of bimodules, respectively. These have 
counterparts in the representation theory of vertex operator algebras. There are 
two commonly used constructions for the tensor product of modules over a vertex 
operator algebra. One of them, developed in the physics literature \cite{nahm26,gaKa},
considers, for a given pair of modules $M$ and $N$, two actions of the vertex 
operator algebra on the vector space \(M\,{\otimes_{\mathbb{C}}}\, N\) and performs a
coequalizer construction resembling the one in the definition of \(\ota\). The other
one \cite{HLZ}, which is e.g.\ briefly summarized in Section 2.2 of \cite{alsW},
starts instead with two actions on the vector space
\(\Hom_{\mathbb{C}}(M\,{\otimes_{\mathbb{C}}}\, N,\mathbb{C})\) and restricts to 
the subspace on which they coincide, similarly as in the equalizer construction
of \(\otimes^A\), and afterwards takes the dual of this subspace.
We refer to \cite{kaRi} for a more detailed account of the two approaches. We expect 
that in favorable circumstances the two constructions can be related through a 
Grothendieck-Verdier structure in a similar way as the tensor products \(\ota\) 
and \(\otimes^A\) of bimodules; we intend to come back to this issue elsewhere.
\end{remark}

Recall from Definition \ref{def:icoHom} that, for $(\calc,\otimes,1)$ a linear
monoidal category and $x,z\iN\calc$, the internal coHom $\icoHomr(z,x) \iN \calc$ 
is the object characterized by the functorial isomorphisms \eqref{eq:def:icoHom}.
As seen in \eqref{eq:iHom=}, for any GV-category internal Homs for $\tpr$ exist, 
and so do internal coHoms for $\tpl\,$. 
In terms of the two tensor products, they can be written as
  \be
  \icoHomr(x,y) = y\tpR G^{-1}x \quad \text{ and }\quad 
  \iHomr(x,y) = y\tpl Gx \,,
  \label{eq:iHom=otiG-1}
  \ee
and as
  \be
  \icoHoml(x,y) = Gx\tpR y \quad \text{ and }\quad 
  \iHoml(x,y) = G^{-1}x \tpl y \,,
  \ee
generalizing the expressions \eqref{eq:iHom/icoHom.as.tpr/tpl}
and \eqref{eq:iHoml/icoHoml.as.tpr/tpl}.

It is an instructive exercise to derive the concrete formulas we gave for the 
monoidal category $A\bimod$ for the tensor product $\otimes^A$ and the internal 
(co)Hom from these general statements. Also, in analogy to \eqref{eq:def:imu} it 
follows immediately that for any object $c\iN\calc$ the internal coHom $\icoHom(c,c)$ 
is a counital coassociative coalgebra for the tensor product $\tpl\,$. By duality, 
this tensor product comes with coevaluations.

 \medskip

As usual, the trivialization of double duals is interesting:

\begin{definition}
A \emph{pivotal} GV-category is a GV-category together with natural isomorphisms
  \be
  \psi_{x,y} \Colon \Hom(x\tpR y,K) \rarr\cong \Hom(y\tpR x, K)
  \ee
such that
  \be
  \psi_{x,y}\circ\psi_{y,x} = \id \qquad \text{and} \qquad
  \psi_{x\tpr y,z}\circ \psi_{y\tpr z,x} \circ \psi_{z\tpr x,y} = \id \,,
  \label{eq:pivax} 
  \ee
or, more explicitly when not suppressing the associatior,
  \be
  \psi_{x\tpr y,z}\circ \alpha_{x,y,z} \circ\psi_{y\tpr z,x} \circ \alpha_{y,z,x}
  \circ \psi_{z\tpr x,y} \circ \alpha_{z,x,y} = \id \,.
  \ee
\end{definition}

We have

\begin{proposition} \cite[Prop.\,6.7]{boDr}
Pivotal structures on $\calc$ are in bijection
to monoidal isomorphisms of functors $\pi\colon \id\To G^2$
whose component $\pi_K \colon K \,{\Rarr\cong}\, G^2(K)$ is given by 
the canonical isomorphism $K \Rarr\cong G 1 \eq G^2G ^{-1} 1 \Rarr\cong G^2 K$.
\end{proposition}

There also exists a notion of a \emph{ribbon} GV-category
\cite[Sect.\,8]{boDr}. It has been shown that ribbon GV-categories 
lead to \emph{ansular functors} \cite{muWo3}, and non-degeneracy conditions 
on the braiding have been identified which guarantee that they lead to
modular functors \cite{brWo}. It should also be appreciated that the representation 
category of a vertex operator algebra to which the HLZ tensor product theory applies 
has a natural structure of a ribbon GV-category \cite{alsW}.
(This fits with the facts that a vertex operator algebra and
its gradewise dual are not necessarily isomorphic as modules, and 
that the tensor product of vertex operator algebra modules need not be exact.)


\subsection{A symmetric formulation of the Eilenberg-Watts theorem}

To connect the abstract considerations of Sections 
\ref{sec:1stdef}\,--\,\ref{sec:2ndtp} with the particular case of $A\bimod$ that
was studied in Section \ref{sec:bimod},
let us focus on the case of a GV-category that is in addition abelian.

\begin{lemma} 
Let $\calc$ be an abelian category that has the structure of a GV-category with 
biadditive tensor product bifunctor $\otimes$. Then $\otimes$ is right exact and 
the second tensor product $\tpl$ is left exact.
\end{lemma}

This statement applies in particular to the abelian category $A\bimod$ 
of finite-dimensional bimodules over finite-dimensional algebras which we considered
in Section \ref{sec:bimod}. In this case we observe, as a by-product, the 
following application of dualities: the GV-structure 
on $A\bimod$ allows one to rewrite the classical Eilenberg-Watts theorem in a form
that treats left exact and right exact functors in a completely symmetric manner.

The classical Eilenberg-Watts theorem describes right exact functors between
categories of modules: Let $\moD A$ and $\moD A'$ be the categories of 
finite-dimensional right modules over finite-dimensional algebras $A$ and $A'$. Then
any right exact functor 
  \be
  F \Colon \moD A \rarr~ \moD A'
  \ee
that commutes with direct sums obeys
  \be
  F(-) \cong {-} \,{\otimes_{\!A}}\, F(A) \,.
  \label{eq:EWright}
  \ee
Here it should be appreciated that $F(A)$ is by definition an $A'$-right module,
while it is also a left $A$-mo\-du\-le by applying $F$ to the endomorphisms of $A$, as
a right module over itself, given by the action of $F$ on left multiplication. A
similar formula is valid for a functor $F'\colon A\Mod\To A'\Mod$ of categories of 
left modules:
  \be
  F'(-) \cong F(A) \,{\otimes_{\!A}}\, {-} \,.
  \ee
  
There is also a version of the Eilenberg-Watts theorem for left exact functors,
in which the role of the tensor product functor in \eqref{eq:EWright} is taken over 
by a Hom functor \cite[Thms.\,2.6\,\&\,2.7]{ivanS}: If (and only if) $H$ is 
left exact, then
  \be
  H \cong \Hom_A(M,-) \quad \text{with} \quad M \eq { (H( {}_{A\!}^{}A_A^*))}^*_{} .
  \ee
(In particular, $H$ is representable.)

Now consider the particular instance
  \be
  G_{\!A,\ko} \Colon A\Mod \rarr~ \moD A
  \ee
of the duality functor given in \eqref{eq:def:GAB} and the analogous functor 
$G_{\ko,A}$ mapping right modules to left modules. Suppose that 
  \be
  H \Colon A\Mod \to A'\Mod
  \ee
is a left exact functor of left modules. Then
  \be
  F := G_{\!A',\ko} \cir H \cir G_{\ko,A} \Colon \moD A \rarr~ \moD A' 
  \ee
is a right exact functor of right modules. 
The Eilenberg-Watts equivalence \eqref{eq:EWright} for this functor $F$ reads
  \be
  G_{\!A',\ko} \cir H \cir G_{\ko,A}(-) 
  \,\cong\, - \,{\otimes_{\!A}} \big( G_{\!A',\ko}\, H\, G_{\ko,A}(A) \big) \,.
  \label{eq:EW4GHG}
  \ee
Now for $M\iN A\Mod$ set $\widetilde M \,{:=}\, G_{A,\ko}(M)\iN \moD A$. Applying
the isomorphism \eqref{eq:EW4GHG} to $\widetilde M$ and then $G_{\ko,A'}$ to the 
so obtained isomorphism gives
  \be
  \begin{aligned}
  H(M) &
  = G_{\ko,A'} \big( \widetilde M \otimes_{\!A} G_{\!A',\ko} H G_{\ko,A}(A) \big) 
  \\[3pt]
  & = G_{\ko,A'} \big( G_{A,\ko}(M) \,{\otimes_{\!A}} G_{\!A',\ko}\, H(A^*) \big)
  = H(A^*) \,{\otimes^{\!A}} M \,.
  \end{aligned}
  \ee
(A corresponding result for comodules is shown in \cite[12.1]{BRwi}.)
This description is beautiful, even for very classical algebra, and promising. In
short, with the help of the duality functor, \emph{left exact} functors can be 
expressed in terms of the second tensor product $\otimes^A$ in completely the 
same way as in the classical Eilenberg-Watts description \eqref{eq:EWright}
of \emph{right exact} functors: apply the functor to the monoidal unit -- this 
time $A^*$ -- and then take the tensor product -- this time $\otimes^A$.

 \smallskip

For \(\varphi\colon A\Rarr{} A'\) a homomorphism of associative \(\ko\)-algebras,
the pullback functor $\mathrm{Res}_\varphi\colon \moD 
        $\linebreak[0]$
        A'\Rarr~ \moD A$ assigns 
to any right \(A'\)-module \((M',\rho_{M'})\) the right \(A\)-mo\-du\-le 
\((M',\rho_{M'}\cir(\varphi\otik \id_{M'}))\) with $A$-action
given by pulling back the $A'$-action along \(\varphi\).
(The notation $\mathrm{Res}$ indicates that if \(\varphi\) is injective, then 
$\mathrm{Res}_\varphi$ is the restriction of the $A'$-action to a subalgebra \(A\).)
For easier distinction we write \(\varphi^\ast(M')\) for the image of \(M'\) under
\(\varphi^\ast\)  even though it coincides with $M'$ as a vector space.
On morphism spaces \(\varphi^\ast\) is then the first inclusion in the sequence
  \be
  \Hom_{A'}(M',N') \subset \Hom_{A}(\varphi^\ast(M'),\varphi^\ast(N'))
  \subset \Hom_{\ko}(M',N') \,.
  \ee
Thus the functor \(\mathrm{Res}_\varphi\) is faithful (albeit not
necessarily full or essentially surjective) and exact. The symmetric formulation
of the Eilenberg-Watts theorem obtained above yields two expressions
  \be
  M' \,{\otimes_{\!A'}}\, \mathrm{Res}_{\varphi}(A') = \mathrm{Res}_\varphi (M')
  = M' \,{\otimes^{\!A'}}\, \mathrm{Res}_{\varphi}({A'}^\ast)
  \ee
for the pullback functor. Denoting the left and right adjoint of 
\(\mathrm{Res}_\varphi\), which in the subalgebra case are the induction and 
coinduction functors, by ${\mathrm{Ind}_\varphi}$ and ${\mathrm{coInd}_\varphi}$,
respectively, the adjunction \ref{eq:Adj-bimod} thus gives 

\begin{corollary} 
The left and right adjoints of the pullback functor $\mathrm{Res}_\varphi$
can be written symmetrically as
  \be
  {\mathrm{Ind}_\varphi}(M) = M \otA \mathrm{Res}_\varphi {(A'^\ast)}^\ast
  \quad \text{~and~} \quad
  {\mathrm{coInd}_\varphi}(M) = M \oaT \mathrm{Res}_\varphi {(A')}^\ast . 
  \ee
\end{corollary}


\section{Distributors}

In a monoidal category with tensor product $\tpr$, different bracketings of
multiple tensor products are related through the associator $\alph \eq \alph^\tpr$.
(A pentagon identity for $\alpha$ ensures coherence.)	
In a GV-category there is another associator $\alph^\tpl$ for multiple
$\tpl$-tensor products, whose components are directly obtained from those of 
$\alpha$ as images under the duality functor $G$. As a consequence, they are 
isomorphisms and also $\alpha^\bullet$ satisfies a pentagon equation.
In addition there are similar coherent families of morphisms that relate
multiple products involving both $\tpl$ and $\tpr$. These are best understood
when viewed from the perspective of module categories, involving in particular
the notion of a weak module functor.
Recall that a module category is a categorification of the notion of a module
over a ring, see e.g.\ \cite[Def.\,7.1.1]{EGno},

\begin{definition}
Let $\calm$ and $\caln$ be left module categories over a monoidal category $\calc$. 
A \emph{lax module functor} from $\calm$ to $\caln$ is a functor
$G\colon \calm \Rarr~ \caln$ endowed with a natural family of morphisms
$g_{c,m}\colon c \actre G(m) \Rarr~ G(c \actrE m)$ (called the lax \emph{module 
constraint} of $G$) such that the appropriate pentagon diagram is fulfilled.
Analogously, for an \emph{oplax module functor} $F\colon \calm \Rarr~ \caln$ there 
is a corresponding coherent family of morphisms $F(c \actrE m) \Rarr~ c \actre F(m)$.
 \\
A lax or oplax module functor is also called a \emph{weak module functor}.
 \\
A \emph{strong module functor} is a lax (or oplax) module functor for which
the lax (or oplax) module constraints are isomorphisms.
\end{definition}

The notion of internal Homs and coHoms generalizes from monoidal categories 
to module categories in a straightforward way.
For $(\calm,\actre)$ a left module category over a monoidal category $(\calc,\tpr)$, 
one defines internal Homs, similarly as $\iHomr$ in Definition \ref{def:iHom} for 
$\calc$ itself, as objects representing the functors 
$c\,{\mapsto}\, \Hom_\calm(c\actrE m,n)$, giving rise to coherent isomorphisms
  \be
  \Hom_\calm(c\actrE m,n) \cong \Hom_\calc(c,\iHomr(m,n)) 
  \label{eq:def:MiHom}
  \ee
for $c\iN\calc$ and $m,n\iN\calm$. Denote, for any object $m$ in a left 
$(\calc,\tpr)$-mo\-du\-le category $\calm$, by 
$\RL_m^\actre \colon {}_{\calc}\calc \Rarr~ \calm$, with ${}_{\calc}\calc$ standing
for the category $\calc$ regarded as a left module over itself,
the strong $(\calc,\tpr)$-module functor obtained by acting on $m$, i.e.
  \be
  \RL_m^\actre(c) := c \actrE m
  \label{eq:def:Lm-actre}
  \ee
for $c\iN\calc$. Then the \emph{internal Hom} $\iHomr(m,-)$ is the right
adjoint of the module functor $\RL_m^\actre$.  Indeed, all all module functors 
$\calc\Rarr~\calm$ are of this type:

\begin{lemma} \label{lemma:generic-module-fun}
Let $\calc$ be a monoidal category, $\calm$ a left $\calc$-module, and
$F \colon \calc \Rarr~ \calm$ a strong module functor. Then 
there exists an object $m \iN \calm$ with a module natural isomorphism
$F\Tocong \RL_{m}^{\actre}$. The object $m$ is unique up to unique isomorphism. 
\end{lemma}

\begin{proof} We set $m \,{:=}\, F(\one)$. From the module constraint $f$ of
$F$ we obtain a natural isomorphism 
$f_{c,1} \colon F(c) \eq F(c \oti \one) \Tocong c \actre F(\one) \eq c \actre m$
for $c \iN \calc$. From the pentagon axiom of the module constraint it follows
that this family of isomorphisms constitutes a module natural isomorphism 
$F\Tocong \RL_m^{\actre}$. By a Yoneda argument, $m$ is unique up to unique isomorphism.
\end{proof}

If $(\calc,\tpr)$ is a GV-category, then the regular left module 
${}_{\calc}\calc$ does admit internal Homs (which are given by the formulas 
\eqref{eq:iHom=}). In contrast, for an arbitrary module category over a 
GV-category, the existence of such right adjoints is, of course, not guaranteed. 
Accordingly we give

\begin{definition}
A \emph{left GV-module category} over a GV-category $\calc$ is a left module 
category $(\calm,\actre)$ over $(\calc,\tpr)$ such that both the action functor 
$\RL_m^\actre \colon {}_{\calc}\calc \Rarr~ \calm$ defined by \eqref{eq:def:Lm-actre} 
and the functor $c \actrE {-} \colon \calm \Rarr~ \calm$ admit a right adjoint,
for every $m\iN\calm$ and every $c \iN \calc$, respectively.
\end{definition}

Note that the GV-structure of $\calc$ is in fact not used in this definition. The
separate terminology `GV-module category' is chosen to indicate that an additional 
condition is imposed. 

Since any left adjoint functor preserves all colimits, for $\calm$ a GV-module 
category over a GV-category $\calc$ the functor $\RL_m^\actre$ is right exact for 
every $m\iN\calm$. By definition, for any $c \iN \calc$ the endofunctor 
$\,c\actrE{-}$ has a right adjoint; 
denote this right adjoint by $R_c\colon \calm\To\calm$. We have

\begin{proposition}  
Let $(\calm,\actre)$ be a left GV-module category over a GV-category $(\calc,\tpr)$. 
Then the bifunctor $\actle \colon \calc\Times\calm\Rarr~\calm$ given by
  \be
  c\actlE m := R_{Gc}(m)
  \label{eq:actlE=RGc}
  \ee
is left exact in each variable and defines a
left module category structure over $(\calc,\tpl)$. 
\end{proposition}  

\begin{proof}
In the defining adjunction $\Hom_\calm(c\actrE n,m) \Cong \Hom_\calm(n,R_c(m))$
of $R_c$, the left hand side is manifestly left exact in $m\iN\calm$, and thus so 
is $\actle$. Left exactness in $c\iN\calc$ follows by right exactness of $\actre$
from the isomorphism
  \be
  \Hom_\calm(m,c\actlE n) = \Hom_\calm(m,R_{Gc}(n)) \cong \Hom_\calm(Gc\actrE m,n) \,.
  \ee
Further, for $m,n\iN\calm$ and $b,c\iN\calc$ there are natural isomorphisms
  \be
  \begin{aligned}
  \Hom_\calm(m,b\actlE (c\actlE n))
  & \cong \Hom_\calm(Gb\actrE m,c\actlE n)
  \Nxl1 &
  \cong \Hom_\calm(Gc\actrE (Gb\actrE m), n)
  \Nxl1 &
  \cong \Hom_\calm((Gc\tpR Gb)\actrE m), n)
  \Nxl1 &
  \cong \Hom_\calm( m,(b\tpL c)\actlE n) \,.
  \end{aligned}
  \ee
{}From these isomorphisms the module constraints for $\actle$ follow by the Yoneda 
lemma. Moreover, the Yoneda embedding transports the pentagon identity for the module
constraint of $\actre$ to the one for the module constraint of $\actle$.
\end{proof}

Setting
  \be
  \RL_m^\actle(c) := c \actlE m
  \ee
we thus obtain for any left GV-module category $(\calm,\actre)$ over a
GV-category $(\calc,\tpr)$ also a $(\calc,\tpl)$-module functor 
$\RL_m^\actle \colon {}_{\calc}^\tpl\calc \Rarr~ \calm$. Since $\RL_m^\actle$ is defined
via the right adjoint of the $(\calm,\actre)$-module functor $\RL_m^\actre$, it
has a left adjoint. Indeed we have natural isomorphisms
  \be
 \begin{aligned}
  \Hom_\calm(n,\RL_{m}^{\actle}(c)) & = \Hom_{\calm}(n, c \actlE m)
  \cong \Hom_\calm(Gc\actrE n, m)   
  \nxl1 &
  \cong \Hom_\calc(Gc, \iHomr(n,m))
  \cong \Hom_\calc(G^{-1}(\iHomr(n,m)),c) \,. 
  \end{aligned}
  \label{eq:cohom-adj}
  \ee
We call the left adjoint of the $(\calc,\tpl)$-module functor $\RL_m^\actle$ the
\emph{internal coHom} and denote it by $\icoHomr(m,-)$.
By definition, there are thus coherent isomorphisms
  \be
  \Hom_\calm(n,c\actlE m) \cong \Hom_\calc(\icoHomr(m,n),c)
  \label{eq:iso:icoHom}
  \ee
for all $c\iN\calc$ and all $m,n\iN \calm$ and, by \eqref{eq:cohom-adj}, 
  \be
  \icoHomr(m,n) \cong G^{-1}(\iHomr(n,m))
  \ee
for all $m,n\iN \calm$. 

It should be appreciated that the adjoint of a strong module functor is, in general, 
only a weak module functor. To understand this, it is convenient to consider 
linear categories and work with module profunctors.
Following \cite[Def.\,2.1]{shimi20} we give
  
\begin{definition}
Let $\calm$ and $\caln$ be left modules over a \ko-linear monoidal category
$\calc$. A \emph{$\calc$-module profunctor} from $\calm$ to $\caln$ is a bilinear
functor $H\colon \calm\opp \Times \caln \,{\rarr~}
	$\linebreak[0]$
\vect_\ko$ together with a family
  \be
  \theta_{m,n,c}\Colon H(m,n) \rarr~ H(c \actrE m, c \actrE n)
  \ee
of morphisms that is natural in $m \iN \calm$ and in $n \iN \caln$, is
dinatural in $c \iN \calc$, and is coherent with respect to
the monoidal structure and the unit of $\calc$, i.e.\ satisfies
$\theta_{m,n,c\otimes d} \eq \theta_{d\actre m,d\actre n,c} \cir \theta_{m,n,d}$ and
$\theta_{m,n,\TI} \eq \id_{H(m,n)}$ for all $m\iN\calm$, $n\iN\caln$ and $c,d\iN\calc$.
(Here, analogously as we do for monoidal categories, we suppress the mixed 
associators and unitors of a module category.)
 \\
A \emph{morphism} $\varphi\colon H_1 \Rarr~ H_2$ \emph{of module profunctors}
between $\calc$-module profunctors $H_1$ and $H_2$ is a natural transformation
that commutes with the respective $($di$)$natural transformations.
\end{definition}

The aspect of module profunctors that is relevant to us is

\begin{lemma} \cite[Lemma\,2.3]{shimi20} \label{lem:transport} 
Let $\calm$ and $\caln$ be left modules over a \ko-linear monoidal category $\calc$, 
and let $F\colon \calm \Rarr~ \caln$ and $G\colon \caln \Rarr~ \calm$ be \ko-linear
functors. There is a canonical bijection between oplax module functor structures 
on $F$ and $\calc$-module profunctor structures on
  \be
  \Hom_{\caln}(F(-),-) \Colon \calm\opp \Times \caln \rarr~ \vect_\ko \,,
  \ee
as well as a canonical bijection between 
lax module functor structures on $G$ and $\calc$-mo\-dule profunctor structures on
  \be
  \Hom_{\calm}(-,G(-)) \Colon \calm\opp \Times \caln \rarr~ \vect_\ko \,.
  \ee
\end{lemma}

This allows one to transport 
lax to oplax module functor structures and vice versa:
     
\begin{corollary} \cite[Lemmas\,2.4\,\&\,2.5]{shimi20} \label{Corollary:adj-lax}
Let $\calm$ and $\caln$ be left modules over a \ko-li\-near monoidal category
$\calc$, and let $F\colon \calm \Rarr~ \caln$ be a linear functor with right adjoint
$G\colon \caln \Rarr~ \calm$. There is a canonical bijection between oplax
$\calc$-module functor structures on $F$ and lax $\calc$-module structures
on $G$, such that the adjunction $\varphi$ with components
  \be
  \varphi_{m,n} \Colon \Hom(F(m),n) \rarr~\Hom(m, G(n))
  \ee
is an isomorphism of $\calc$-module profunctors. If $F$ is even a strong module 
functor, then $G$ has a unique structure of a lax $\calc$-module functor such that 
the unit and counit of the adjunction $\varphi$ are module natural transformations. 
\end{corollary}

It follows in particular, by considering the module functor $F \eq \RL_{m}^{\actre}$,
that for every object $m$ in a left $\calc$-module category $\calm$ the internal
Hom $\iHomr(m,-)$ is a lax module functor. Hence there are natural coherent morphisms
  \be
  \deltar cmn\Colon c \tpr \iHomr(m,n) \rarr~ \iHomr(m,c \actrE n)
  \ee
for $c\iN\calc$ and $m,n\iN\calm$.

\begin{remark} \label{rem:BOur3takeu3}
The observation that internal Homs are not necessarily strong module functors is,
of course, not new, e.g.\ it is used in the definition of rigidity in
\cite[Def.\,9.1.2]{GAroz}. Conditions for being strong have been considered in the
literature. For instance, for finitely generated projective modules 
the internal Hom is a strong module functor \cite[p.\,269]{BOur3}.
For internal coHoms, which naturally appear in the study of categories of comodules
over a coalgebra, a similar statement, namely that one has strong module functors for 
quasi-finite and injective comodules, has been obtained in \cite[Prop.\,1.14]{takeu3}.
See also Lemma \ref{lem:iHom-strong} for equivalent conditions for internal Homs 
and coHoms in categories of bimodules to be strong module functors.
\end{remark}

The particular case of interest to us is that $\calc$ is a GV-category and that
$\calm \eq {}_{\calc}\calc$ is the regular left $\calc$-module,
for which according to \eqref{eq:iHom=otiG-1} we have $\iHomr(x,y) \eq y \tpl G(x)$.
In this case we obtain

\begin{lemma}
Let $\calc$ be a \ko-linear GV-category. Then there is a family
  \be
  \deltar cxy :~~ c \tpR \iHomr(x,y) = c \tpr (y \tpl G(x))
  \rarr~ (c \tpR y ) \tpl G(x) = \iHomr(x,c \tpR y)
  \label{eq:def:deltal}
  \ee
of morphisms, for $c,x,y\iN\calc$, which endows 
the internal Hom $\iHomr(x,y) \eq y \tpl G(x)$
with a lax module functor structure for $\calc$ as the regular left $\calc$-module.
\end{lemma}

While $\iHomr(x,-) \eq {-} \tpL G(x)$ is, in general, only a weak module functor for 
the right exact tensor product $\tpr$, it is a strong module functor for the left
exact tensor product $\tpl$. Indeed, the module functor
  \be
  I_x := \iHomr(x,-) \Colon {}_{\calc}^{\tpl}\calc \rarr~ {}_{\calc}^{\tpl}\calc \,,
  \ee
with ${}_{\calc}^{\tpl}\calc$ the regular $(\calc,\tpl)$-left module, is strong 
because the associator $\alph^\tpl$ provides us with isomorphisms
  \be
  I_y(c \tpL x) = (c\tpL x) \tpL G(y) \rarr{\cong} c \tpL (x \tpL G(y))
  = c \tpL I_{y}(x) \,.
  \ee
Invoking Lemma \ref{lem:transport}, this implies that the left adjoint 
$\RL_y \eq {-}\tpR y\colon {}_{\calc}^{\tpl}\calc \Rarr~ {}_{\calc}^{\tpl}\calc$ 
of $I_y$ is an oplax module functor. Accordingly there are coherent morphisms
  \be
  \begin{aligned}
  \deltal yxc \Colon (G\inv(x) \tpL y) \tpR c =~ & \RL_c(G\inv(x) \tpL y)
  \\
  & \rarr{\phantom{ii}} G\inv(x) \tpL \RL_c(y) = G\inv(x) \tpL (y \tpR c) \,.
  \end{aligned}
  \label{eq:def:deltar}
  \ee

Borrowing terminology from the theory of linearly distributive categories
(see e.g.\ \cite{coSe3,past,haLem}), we give

\begin{definition} \label{def:distributors}
The natural transformations $\deltaL$ and $\deltaR$ introduced in 
\eqref{eq:def:deltal} and \eqref{eq:def:deltar} are called the left and right 
\emph{distributors} of the GV-category $\calc$, respectively.
\end{definition}

In fact, as stated in the literature \cite{past,haLem}, GV-categories are the
same as linearly distributive categories with negation. The interpretation of
distributors in terms of GV-structures has, however, (to the best of our knowledge)
not been established explicitly. In our approach, the definition of distributors 
via weak module functors implies

\begin{proposition}
The distributors $\deltaL$ and $\deltaR$ 
satisfy all the compatibility conditions with the unitors for $\TI^\tpr$ and
$\TI^\tpl$ $($four mixed triangle identities$)$ and with the associators
$\alph^\tpr$ and $\alph^\tpl$ $($four mixed pentagon identities$)$
that the distributors in a linearly distributive category have to obey.
\end{proposition}

\begin{proof}
Since by definition a weak module functor satisfies an appropriate pentagon 
diagram, the two paths in the pentagon diagram
that is built from the associator $\alpha^\tpr$ and the weak module structure 
\ref{eq:def:deltal}, which induces the left distributor, commute:
  \be
  \begin{aligned}
  (c \tpR d) \tpR\iHomr(x,y) & = (c \tpR  d) \tpR  (y \tpL G(x))
  \nxl1
  & \!\!\! \begin{tikzcd}[column sep=1.5em] {} \ar[r,yshift=2pt,bend left=4]{}
     \ar[r,yshift=-3pt,bend right=4]{} & {} \end{tikzcd} \!\!\!
     (c \tpR (d \tpR y)) \tpl G(x) = \iHomr(x,c \tpR (d \tpR y)) \,.
  \end{aligned}
  \label{eq:pentag1}
  \ee
Further, for all $x_{1},x_{2} \iN \calc$ there is a canonical isomorphism 
  \be
  L^{}_{x_2} \circ L^{}_{x_1} \cong L^{}_{x_1 \tpr x_2}
  \ee
of strong module functors. This isomorphism induces an isomorphism of the respective 
right adjoint weak module functors, and thus translates into a commuting diagram
  \be
  c \tpR ((y \tpL G(x_2)) \tpL G(x_1))
  \!\!\! \begin{tikzcd}[column sep=1.5em] {} \ar[r,yshift=2pt,bend left=4]{}
     \ar[r,yshift=-3pt,bend right=4]{} & {} \end{tikzcd} \!\!\!
  ( c \tpR y) \tpL (G(x_2) \tpL G(x_1)) \,.
  \label{eq:pentag2}
  \ee
Two further pentagon diagrams, which are mirror images of \eqref{eq:pentag1} and
\eqref{eq:pentag2} and involve the right distributor commute by the analogous 
arguments for the strong module functors $L^{\tpl}_{x}$. 
Inspection shows that these pentagons are the same as those obeyed by the 
distributors in a linearly distributive category, as given e.g.\ in
\cite[Sect.\,2.1.3]{coSe3}. The triangle diagrams involving the unitors 
(see e.g.\ \cite[Sect.\,2.1.2]{coSe3})
are immediate.
\end{proof}
    
It can also be shown that the two additional pentagon identities that are valid in a
linearly distributive category, each of which involves both $\deltaL$ and $\deltaR$,
are fulfilled as well. However, the proof of this statement we know of is less
conceptual and considerably more indirect; we refrain from presenting it here.
Also note that, as a consequence of their construction via \emph{weak} module 
functors, the distributors $\deltaR$ and $\deltaL$ are, in general, not
isomorphisms. They \emph{are} isomorphisms if and only if the category is rigid,
see Lemma \ref{lem:iHom-strong}.


\section{Subcategories of rigid objects}

Next we characterize, for right closed monoidal categories $\calc$, those objects
$x \iN \calc$ for which $\iHomr(x,-)$ is a \emph{strong} module functor, i.e.\ for
which the coherence morphisms ${}^{\ihl}\delta^x_{y,z}$
in \eqref{eq:def:deltal} are \emph{iso}morphisms for all objects $y,z \iN \calc$.
We start with the standard observation (compare e.g.\ \cite[Prop.\,2.10.8]{EGno})
that the monoidal equivalence of $\calc$ and
$\End_{\calc}({}_{\calc}\calc)$ directly implies

\begin{lemma} \label{Lemma:Duals-adjoint}
Let $(\calc,\tpr)$ be a monoidal category and let $x\iN\calc$. An object
$x^{\vee}\iN\calc$ is a right dual of $x$ if and only if the functor 
${-} \tpr x^{\vee}$ is right adjoint to the functor ${-} \tpr x$ as a module functor.
\end{lemma}

Note that it does not suffice to merely require that ${-} \tpr x^{\vee}$ is right 
adjoint to ${-} \tpr x$ as a linear functor (a counter example is given in \cite{haZo2}). 
A statement analogous to Lemma \ref{Lemma:Duals-adjoint} holds for left duals.

\begin{proposition} \label{Proposition:Strong-duals}
Let $(\calc,\tpr)$ be a right closed monoidal category. For $x \iN \calc$ the lax 
module functor $\iHomr(x,-) \colon \calc \Rarr~ \calc$ is a strong module functor if 
and only if $x$ has a right dual object $x^\vee\!$. Moreover, if this is the case, 
then $x^{\vee}\,{\cong}\, \iHomr(x,\TI)$ as objects, and 
  \be
  \iHomr(x,-) \cong {-} \tpr x^{\vee}
  \ee
as module functors.
\end{proposition}

\begin{proof}
Assume that the object $x$ has a right dual $x^\vee$. Then the left 
$\calc$-mo\-du\-le functor 
$\RL_x\colon \calc \Rarr~ \calc$ with $\RL_x(y) \eq y \tpR x$ has as a right adjoint
the left $\calc$-mo\-du\-le functor $\RL_{x^{\vee}}\colon \calc \Rarr~ \calc$, in such
a way that the unit and counit of the adjunction
$\Hom_{\calc}(\RL_x(y),z)\,{\cong}\, \Hom_{\calc}(y,\RL_{x^\vee}(z))$
consist of module natural transformations. Thus it follows from
Corollary \ref{Corollary:adj-lax} that $\iHomr(x,-) \eq {-} \tpR x^\vee$ as 
lax module functors. This, in turn, implies that $\iHomr(x,-)$ is in fact a 
strong module functor, with the associator of $\calc$ as module constraint.
 \\
Conversely, assume that $\iHomr(x,-)\colon \calc \Rarr~ \calc$ is a strong module
functor. By Lemma \ref{lemma:generic-module-fun}
it then follows that there is a natural isomorphism
  \be
  \iHomr(x,-) \cong {-} \tpR \iHomr(x,\TI)
  \ee
of module functors. Thus ${-} \tpR \iHomr(x,\TI)$ is right adjoint to ${-} \tpR x$
as a module functor, and so by Lemma \ref{Lemma:Duals-adjoint} the object
$\iHomr(x,\TI) \iN \calc$ is a right dual of $x$.
\end{proof}

Now the class of objects of $(\calc,\tpr)$ that admit a right dual is closed under 
the monoidal product $\tpr$. Thus we have

\begin{corollary} \label{Cor:Mon-subcat-dualizable}
The full subcategory on those objects $y \iN \calc$ for which $\iHomr(y,-)$ is a
strong module functor is a unital monoidal subcategory of $\calc$.
\end{corollary}

Note that, as for instance the category $\AA_3\bimod$ of bimodules over the algebra
\eqref{eq:def:A3} illustrates, for a GV-category the dualizing object $K$ need not
be contained in this subcategory. Also, clearly, there are analogues for left
closed monoidal categories and for internal coHoms:

\begin{lemma}
Let $(\calc,\tpr)$ be a monoidal category.
 \begin{enumerate} \addtolength{\itemsep}{1.8pt}
 \item
Assume that the monoidal category $(\calc,\tpr\opp)$ obtained by reversing the
monoidal structure is left closed. Then for every $x \iN \calc$ the left internal
Hom $\iHoml(x,-)$ is a lax right module functor, with module constraint
$\iHoml(x, y ) \tpR z \Rarr~ \iHoml(x,y \tpR z )$.
 \\
Moreover, $\iHoml(x,-)$ is strong if and only if $x$ has a left dual.
 \item
Assume that the opposite category $(\calc\opp_{},\tpr)$ is left closed. A right
internal coHom of $(\calc,\tpr)$ is a right internal Hom of 
$(\calc\opp_{},\tpr)$.
As a consequence, $\icoHomr(x,-)$ is an oplax left module functor with 
module constraint $\icoHomr(x,y \oti z) \Rarr~ y \otimes \icoHomr(x,z)$.
 \\
Moreover, $\icoHomr(x,-)$ is strong if and only if $x$ has a left dual.
 \item
Assume that the monoidal category $(\calc\opp_{},\tpr\opp)$ is left closed. Then
the left internal coHom $\icoHoml$ is an oplax right module functor with module
constraint $\icoHoml(x,y \tpR z) \rarr~ \icoHoml(x,y) \tpR z$.
 \\
Moreover, $\icoHoml(x,-)$ is strong if and only if $x$ has a right dual.
 \end{enumerate}
\end{lemma}

As an illustration, for categories of bimodules we get

\begin{lemma} \label{lem:iHom-strong}
Let $A$ be a finite-dimensional \ko-algebra, and let ${}_{A}M_{A} \iN A\bimod$
be a finite-dimensional $A$-bimodule.
 \\[2pt]
The following statements are equivalent:
  \begin{itemize} \addtolength{\itemsep}{1.5pt}
  \item[(i)]
$\iHomr(M,-)$ is a strong module functor.
  \item[(ii)]
$M$ has an $\ota$-right dual.
  \item[(iii)]
$M_{A}$ is  projective as a right $A$-module.
  \item[(iv)]
${}_{A}(M^\ast_{})$ is injective as a left $A$-module.
  \item[(v)]
For all $X,Y\iN A\bimod$ the distributor
  \nxl1
$\deltar XMY \colon
X \ota (Y {\otimes^A} M^\ast_{}) \Rarr~ (X \ota Y) \,{\otimes^A} M^\ast_{}$
is an isomorphism.
  \end{itemize}
Likewise, the following statements are equivalent:
  \begin{itemize} \addtolength{\itemsep}{1.5pt}
  \item[(i$'$)]
$\iHoml(M,-)$ is a strong module functor.
  \item[(ii$'$)]
$M$ has an $\ota$-left dual.
  \item[(iii$'$)]
${}_{A}M$ is projective as a left $A$-module.
  \item[(iv$'$)]
$(M^\ast_{})_A$ is injective as a right $A$-module.
  \item[(v$'$)]
For all $X,Y\iN A\bimod$ the distributor
  \nxl1
$\deltal XMY \colon
(M^\ast_{} {\otimes^A} X) \,{\ota} Y \Rarr~ M^\ast_{} {\otimes^A} (X \,{\ota} Y)$
is an isomorphism.
  \end{itemize}
\end{lemma}

\begin{proof}
(ii) follows from (i) with Proposition \ref{Proposition:Strong-duals}. The 
equivalence (ii)$\,\Longleftrightarrow\,$(iii) is standard (compare e.g.\ 
\cite{ncat}).
Since $(-)^{*}$ is an antiequivalence, (iv) is equivalent to (iii). Finally,
equivalence of (v) and (i) follows again from Proposition 
\ref{Proposition:Strong-duals}.
\end{proof}

We can also determine the duals explicitly. For instance, the $\ota$-right dual 
of ${}_{A}M_{A}$ is
  \be
  G^{-1}(M \tpR K) = (M \otA A^\ast)^\ast_{} = \Hom_{A}(M_{A},A_{A}) \,,
  \label{eq:right-dual-M}
  \ee
in accordance with \eqref{eq:iHomrNX}.


\section{Distributors for bimodules}

Let us now describe in detail the distributors for the case of categories of
finite-dimensional bimodules over finite-dimensional $\ko$-algebras considered 
in Section \ref{sec:bimod}. Recall that we assume
all algebras $A$, $B$, etc.\ as well as all bimodules to be finite-dimensional.

Making use of equalizer inclusions and coequalizer surjections,
for any triple \(X,Y,Z\) of \((A,A)\)-bimodules we will define four
\((A,A)\)-bimodule homomorphisms
  \be
  \begin{aligned}
  & \partial^{\ihl}_{X,Y,Z} \,,\, \widetilde\partial^{\ihl}_{X,Y,Z} \Colon
  X \otA (Y \oaT Z) \rarr~ (X \otA Y) \oaT Z
  \\
  \text{and}\qquad & \partial^{\ihr}_{X,Y,Z} \,,\, \widetilde\partial^{\ihr}_{X,Y,Z} \Colon
  (X\oaT Y) \otA Z \rarr~ X \oaT (Y \otA Z)
  \end{aligned}
  \label{eq:leftrightdistributors}
  \ee
which, up to vector space dualities, are natural candidates for the distributors
introduced in Definition \ref{def:distributors}.
A priori, these are all different, but we will show that in fact
$\widetilde\partial^{\ihl}_{X,Y,Z} \eq \partial^{\ihl}_{X,Y,Z}$ and
$\widetilde\partial^{\ihr}_{X,Y,Z} \eq \partial^{\ihr}_{X,Y,Z}$.

Let us first present the definition of $\partial^{\ihl}_{X,Y,Z}$ in some detail.
Consider the diagram
  \be
  \begin{tikzcd}[column sep=3.1em,row sep=2.2em]
  X\otik(Y\otik Z) \ar{r}{\alpha\inv_{\text{vect}}}[swap]{\cong} & (X\otik Y)\otik Z
  \\
  X\otik(Y\oaT Z) \ar[ur, bend right=11]
  \ar[u,hookrightarrow,"X \otimes_\ko \imath^{}_{Y,Z}\!"]
  & (X\otik Y)\oaT Z \ar[hookrightarrow]{u}[swap]{\imath^{\ihl}_{X,Y,Z}}
  \end{tikzcd}
  \ee
where we explicitly display the associator $\alpha_{\text{vect}}^{}$ of the
category of finite-dimensional vector spaces, and where
  \be
  \begin{aligned}
  & \imath^{}_{Y,Z} \Colon Y \oaT Z \hookrightarrow Y \otik Z
  \\
  \text{and}\qquad & 
  \imath^{\ihl}_{X,Y,Z} \Colon (X\otik Y) \oaT Z \hookrightarrow (X\otik Y) \otik Z
  \end{aligned}
  \label{eq:notation-12}
  \ee
are the equalizer inclusions of $Y \oaT Z$ into $Y \otik Z$ and of
$(X\otik Y) \oaT Z$ into $(X\otik Y) \otik Z$, respectively.
Note that since the tensor product of vector spaces is exact,
the map $X \otik \imath^{}_{Y,Z}$ is injective. Explicitly we have
  \be
  \begin{aligned}
  \alpha_{\text{vect}}^{-1} \circ (X \otik \imath^{}_{Y,Z}) \Colon
  X\otik (Y\oaT Z) & \rarr~ (X\otik Y)\otik Z \,,
  \\
  \sum_{i} x\oti (y_i\oti z_i) & \xmapsto{~~~} \sum_{i} (x\oti y_i)\oti z_i \,,
  \end{aligned}
  \label{eq:alphainv.Ximath}
  \ee
where $x\iN X$ and we write a generic element of \(Y\oaT Z\) as a finite sum
$\sum_i y_i\otik z_i$ with $y_i\iN Y$ and $z_i\iN Z$, obeying
$\sum_i y_i.a\otik z_i \eq \sum_i y_i\otik a.z_i$ for all $a\iN A$.
The map \eqref{eq:alphainv.Ximath} is cobalanced, i.e.\
$\sum_{i} (x\oti y_i)\oti a.z_i \eq \sum_{i} (x\oti y_i.a)\oti z_i
\eq \sum_{i} (x\oti y_i).a\oti z_i$ for all $a\iN A$; hence by the universal property 
of equalizers there exists a unique homomorphism 
$\gamma^{\ihl}_{X,Y,Z} \colon X\otik(Y\oaT Z)\Rarr~ (X\otik Y)
       $\linebreak$
       {\oat}\, Z$ such that the diagram
  \be
  \begin{tikzcd}[column sep=3.5em,row sep=2.2em]
  X\otik(Y\otik Z) \ar{r}{\alpha\inv_{\text{vect}}}[swap]{\cong} & (X\otik Y)\otik Z
  \\
  X\otik(Y\oaT Z) \ar[r,dashed,"\gamma^{\ihl}_{X,Y,Z}"]
  \ar[u,hookrightarrow,"X \otimes_\ko \imath^{}_{Y,Z}\!"]
  & (X\otik Y)\oaT Z \ar[hookrightarrow]{u}[swap]{\imath^{\ihl}_{X,Y,Z}}
  \end{tikzcd}
  \label{eq:quare1a}
  \ee
commutes. Since $\imath^{\ihl}_{X,Y,Z}$ and $X \otik \imath^{}_{Y,Z}$ are injective 
and $\alpha^{}_{\text{vect}}$ is an isomorphism, and since the images of
$\imath^{}_{Y,Z}$ and of \eqref{eq:alphainv.Ximath} coincide,
$\gamma^{\ihl}_{X,Y,Z}$ is an isomorphism. Next consider the diagram
  \be
  \begin{tikzcd}[column sep=3.4em,row sep=2.2em]
  X\otik(Y\oaT Z) \ar[twoheadrightarrow]{d}[swap]{\pi^{\ihl}_{X,Y,Z}}
  \ar{r}{\gamma^{\ihl}_{X,Y,Z}}[swap]{~~\cong} \ar[dr, bend left=7]{}
  & (X\otik Y)\oaT Z \ar{d}{\pi^{}_{X,Y} \oat Z}
  \\
  X\otA(Y\oaT Z) & (X\otA Y)\oaT Z
  \end{tikzcd}
  \label{eq:ldistlowersquare}
  \ee
where
  \be
  \begin{aligned}
  & \pi^{}_{X,Y} \Colon X \otik Y \twoheadrightarrow X \otA Y
  \\
  \text{and}\qquad &
  \pi^{\ihl}_{X,Y,Z} \Colon X \otik (Y \oaT Z) \twoheadrightarrow X \otA (Y \oaT Z)
  \label{eq:notation-34}
  \end{aligned}
  \ee
are the coequalizer surjections from $X \otik Y$ onto $X \otA Y$ and from
$X \otik (Y \oaT Z)$ onto $X \otA (Y 
       $\linebreak$
       {\oat}\, Z)$, respectively.
Note that since the functor \(-\oaT Z\) need not be right exact, 
$\pi^{}_{X,Y} \oaT Z$ need not be surjective. Explicitly we have
  \be
  \begin{aligned}
  (\pi^{}_{X,Y} \oaT Z) \cir \gamma^{\ihl}_{X,Y,Z} \Colon
  X\otik(Y\oaT Z) & \rarr~ (X\otA Y)\oaT Z \,,
  \nxl4
  \sum_i x\otik(y_i\otik z_i) & \xmapsto{~\phantom{\cong}} \sum_i [x\otik y_i]\otik z_i \,,
  \label{eq:balancedmorformula}
  \end{aligned}
  \ee
where $[x\otik y]$ is the element of $X \otA Y$ represented by $x \otik y$ with 
$x\iN X$ and $y\iN Y$. The map \eqref{eq:balancedmorformula} is balanced, i.e.\ 
$\sum_i [x.a \otik y_i]\otik z_i \eq \sum_i [x \otik a.y_i]\otik z_i$
for all $a\iN A$, By the universal property of the coequalizer defining the
\(\ota\)-tensor product, there is thus a unique map
$\partial^{\ihl}_{X,Y,Z} \colon X\otA(Y\oaT Z) \,{\rarr~}\, (X\otA Y)\oaT Z$
such that the diagram
  \be
  \begin{tikzcd}[column sep=3.9em,row sep=2.2em]
  X\otik(Y\oaT Z) \ar[twoheadrightarrow]{d}[swap]{\pi^{\ihl}_{X,Y,Z}}
  \ar{r}{\gamma^{\ihl}_{X,Y,Z}}[swap]{\cong}
  & (X\otik Y)\oaT Z \ar{d}{\pi^{}_{X,Y} \oat Z}
  \\
  X\otA(Y\oaT Z) \ar[dashed]{r}{\partial^l_{X,Y,Z}} & (X\otA Y)\oaT Z
  \end{tikzcd}
  \label{eq:square1b}
  \ee
commutes. Explicitly we have
  \be
  \partial^{\ihl}_{X,Y,Z} \Colon
  \sum_i [x\otik(y_i\otik z_i)] \xmapsto{\phantom{\cong }}\sum_i [x\otik y_i]\otik z_i \,.
  \label{eq:firstdistformula}
  \ee
The other three maps in \eqref{eq:leftrightdistributors} are constructed similarly,
leading to

\begin{definition}
Let $X$, $Y$ and $Z$ be \((A,A)\)-bimodules.
 \\
The homomorphism $\partial^{\ihl}_{X,Y,Z} \colon X\otA(Y\oaT Z) \Rarr~ (X\otA Y)\oaT Z$
is the one that is determined by the commutativity of the two squares
\eqref{eq:quare1a} and \eqref{eq:square1b}.
 \\
The homomorphisms 
$\widetilde\partial^{\ihl}_{X,Y,Z} \colon X\otA(Y\oaT Z) \Rarr~ (X\otA Y)\oaT Z$ and
$\partial^{\ihr}_{X,Y,Z}, \widetilde\partial^r_{X,Y,Z} \colon (X\oaT Y) \,{\ota} 
        $\linebreak[0]$ 
        Z \Rarr~ X \oaT (Y \otA Z)$
are the ones that are determined by the commutativity of the following
pairs of squares, respectively (using self-explanatory notation similar to the 
one in \eqref{eq:notation-12} and \eqref{eq:notation-34}):
 \\[2pt]
$\bullet~\widetilde\partial^l_{X,Y,Z}$:
  \be
  \begin{aligned}
  \begin{tikzcd}[column sep=3.4em,row sep=2.2em]
  X\otik(Y\otik Z) \ar{r}{\alpha\inv_{\text{vect}}}[swap]{\cong}
  \ar[twoheadrightarrow]{d}[swap]{\widetilde\pi^{\ihl}_{X,Y,Z}}
  & (X\otik Y)\otik Z \ar[twoheadrightarrow]{d}{\pi^{}_{X,Y} \otimes_\ko Z}
  \\
  X\otA(Y\otik Z) \ar[dashed]{r}{\widetilde\gamma^{\ihl}_{X,Y,Z}}[swap]{\cong}
  & (X\otA Y)\otik Z
  \end{tikzcd}
  \nxl2
  \text{and}\qquad
  \begin{tikzcd}[column sep=3.4em,row sep=2.2em]
  X\otA(Y\otik Z) \ar{r}{\widetilde\gamma^{\ihl}_{X,Y,Z}}[swap]{\cong}
  & (X\otA Y)\otik Z
  \\
  X\otA(Y\oaT Z) \ar{u}{X\ota\imath^{}_{Y,Z}}
  \ar[dashed]{r}{\widetilde\partial^{\ihl}_{X,Y,Z}}
  & (X\otA Y)\oaT Z \ar[hookrightarrow]{u}[swap]{\widetilde{\imath}^{\ihl}_{X,Y,Z}}
  \end{tikzcd}
  \end{aligned}
  \label{eq:def:widetildepartial-l}
  \ee
$\bullet~\partial^{\ihr}_{X,Y,Z}$:
  \be
  \begin{aligned}
  \begin{tikzcd}[column sep=3.4em,row sep=2.2em]
  (X\otik Y)\otik Z \ar{r}{\alpha^{}_{\text{vect}}}[swap]{\cong} & X\otik(Y\otik Z)
  \\
  (X\oaT Y)\otik Z \ar[dashed]{r}{\gamma^{\ihr}_{X,Y,Z}}[swap]{\cong}
  \ar[hookrightarrow]{u}{\imath^{}_{X,Y}\otimes_\ko Z}
  & X\oaT(Y\otik Z) \ar[hookrightarrow]{u}[swap]{\imath^{\ihr}_{X,Y,Z}}
  \end{tikzcd}
  \nxl2
  \text{and}\qquad
  \begin{tikzcd}[column sep=3.4em,row sep=2.2em]
  (X\oaT Y)\otik Z \ar{r}{\gamma^{\ihr}_{X,Y,Z}}[swap]{\cong}
  \ar[twoheadrightarrow]{d}[swap]{\pi^{\ihr}_{X,Y,Z}}
  & X\oaT(Y\otik Z) \ar{d}{X \otimes^A \pi^{}_{Y,Z}}
  \\
  (X\oaT Y)\otA Z \ar[dashed]{r}{\partial^{\ihr}_{X,Y,Z}} &  X\oaT(Y\otA Z)
  \end{tikzcd}
  \end{aligned}
  \ee
$\bullet~\widetilde\partial^{\ihr}_{X,Y,Z}$:
  \be
  \begin{aligned}
  \begin{tikzcd}[column sep=3.4em,row sep=2.2em]
  (X\otik Y)\otik Z \ar{r}{\alpha^{}_{\text{vect}}}[swap]{\cong}
  \ar[twoheadrightarrow]{d}[swap]{\widetilde\pi^{\ihr}_{X,Y,Z}}
  & X\otik(Y\otik Z) \ar[twoheadrightarrow]{d}{X\otimes_\ko \pi^{}_{Y,Z}}
  \\
  (X\otik Y)\otA Z \ar[dashed]{r}{\widetilde\gamma^{\ihr}_{X,Y,Z}}[swap]{\cong}
  & X\otik(Y\otA Z)
  \end{tikzcd}
  \nxl2
  \text{and}\qquad
  \begin{tikzcd}[column sep=3.4em,row sep=2.2em]
  (X\otik Y)\otA Z \ar{r}{\widetilde\gamma^{\ihr}_{X,Y,Z}}[swap]{\cong}
  & X\otik(Y\otA Z)
  \\
  (X\oaT Y)\otA Z \ar[dashed]{r}{\widetilde\partial^r_{X,Y,Z}}
  \ar{u}{\imath^{}_{X,Y} \ota Z}
  & X\oaT (Y\otA Z) \ar[hookrightarrow]{u}[swap]{\widetilde\imath^{\ihr}_{X,Y,Z}}
  \end{tikzcd}
  \end{aligned}
  \ee
\end{definition}

The so defined homomorphisms indeed provide us with the distributors:

\begin{proposition}
Let $X$, $Y$ and $Z$ be \((A,A)\)-bimodules and let
  \be
  \begin{aligned}
  \tensor*[^\ihr]{\delta}{^Z_X_,_Y}\Colon X \otA \iHomr(Z,Y) & \rarr~ \iHomr(Z,X\otA Y)
  \nxl2
  \text{and} \qquad
  \tensor*[^\ihl]{\delta}{^Z_X_,_Y}\Colon \iHoml(Z,X) \otA Y & \rarr~ \iHoml(Z,X\otA Y)
  \end{aligned}
  \ee
be the respective lax module functor structures on the functors \(\iHomr(Z,-)\)
and \(\iHoml(Z,-)\), as defined in \eqref{eq:def:deltal} and \eqref{eq:def:deltar}.
Then we have
  \be
  \partial^{\ihr}_{X,Y,Z} = \tensor*[^\ihl]{\delta}{^{X^\ast}_Y_,_Z}
  = \widetilde{\partial}^\ihr_{X,Y,Z}
  \qquad \text{and} \qquad
  \partial^{\ihl}_{X,Y,Z} = \tensor*[^\ihr]{\delta}{^{Z^\ast}_X_,_Y}
  = \widetilde{\partial}^\ihl_{X,Y,Z} \,.
  \ee
Further, $\partial^{\ihr}$ and $\partial^{\ihl}$ are explicitly given by
  \be
  \begin{aligned}
  & \partial^\ihr_{X,Y,Z}\big( \sum_i [(x_i\otik y_i)\otik z]) 
  = \sum_i x_i\oti [y\otik z]
  \\
  \text{and} \qquad
  & \partial^\ihl_{X,Y,Z}\big( \sum_j [x\otik (y_j\otik z_j)])
  = \sum_j [x\otik y_j]\otik z_j \,.
  \end{aligned}
  \label{eq:distformula}
  \ee
\end{proposition}

\begin{proof}
We show that $\partial^\ihl_{X,Y,Z}$, $\widetilde{\partial}^\ihl_{X,Y,Z}$ and 
$\tensor*[^\ihl]{\delta}{^{Z^\ast}_X_,_Y}$ are all given by the second formula in
\eqref{eq:distformula}. We have already obtained that
formula for \(\partial^\ihl_{X,Y,Z}\) in \eqref{eq:firstdistformula}. The same result
is found for $\widetilde{\partial}^\ihl_{X,Y,Z}$ when performing an 
analogous calculation as the one leading to \eqref{eq:firstdistformula} instead for
the diagrams \eqref{eq:def:widetildepartial-l}. To derive the formula also for
$\tensor*[^\ihl]{\delta}{^{Z^\ast}_X_,_Y}$ we invoke
Lemma \ref{lem:transport} and Corollary \ref{Corollary:adj-lax}
to
compute the module functor structure from the associator $\alpha^\otimes$ of the
$\ota$-tensor product. Consider the adjoint functors 
\(\RL_{Z^\ast}(-) \eq {-}\otA Z^\ast\) and \(I_{Z^\ast_{}}(-) \eq \iHomr(Z^\ast,-)\). 
The associator \(\alpha^\otimes\) provides a strong module functor structure on 
$\RL_{Z^\ast}$ and can be used to construct a transformation
  \be
  \begin{aligned}
  \theta^\otimes_{U,Y,X} \Colon \Hom(\RL_{Z^\ast_{}}(U),Y)
  & \rarr~ \Hom(\RL_{Z^\ast_{}}(X\otimes_A U),X\otimes_A Y) \,,
  \\
  \psi &
  \xmapsto{~~~} (\id_C\otA \psi)\cir (\alpha^{\otimes}_{X,U,Z^\ast})^{-1}_{}
  \end{aligned}
  \ee
that is natural in $U$ and $Y$ as well as dinatural in $X$.
 \\
The internal Hom
adjunction \(\phi_{V,W}\colon \Hom(\RL_{Z^\ast}(V),W) \Tocong \Hom(V,I_{Z^\ast}(W))\)
can then be used to construct a (di)na\-tu\-ral transformation
  \be
  \theta^\delta_{U,Y,X} \Colon
  \Hom(U,I_{Z^\ast_{}}(Y))\rarr~ \Hom(X\otA U,I_{Z^\ast_{}}(X\otimes_A Y))
  \ee
by setting \(\theta^\delta_{U,Y,X} \,{:=}\, \phi^{}_{X\otimes_A U,X\otimes_A Y}\cir
\theta^\otimes_{U,Y,X}\cir \phi^{-1}_{U,Y}\). That is,
\(\theta^\delta_{U,Y,X}\) is exactly such that the diagram
  \be
  \begin{tikzcd}[column sep=3.4em,row sep=2.5em]
  \Hom(\RL_{Z^\ast}(U),Y) \ar{r}{\theta^\otimes_{U,Y,X}} \ar{d}[swap]{\phi^{}_{U,V}}
  & \Hom(\RL_{Z^\ast}(X\otimes_A U),X\otA Y) \ar[d,"\phi^{}_{X\ota U,X\ota V}"]
  \\
  \Hom(U,I_{Z^\ast_{}}(Y)) \ar[r,"\theta^\delta_{U,Y,X}"]
  & \Hom(X\otA U,I_{Z^\ast_{}}(X\otA Y))
  \end{tikzcd}
  \ee
commutes. The transported lax module structure on \(I_{Z^\ast}\) is then
  \be
  {}_{}^{\ihl}\delta^{Z^\ast}_{X,Y} = \theta^\delta_{I_{Z^\ast}(Y),Y,X}
  (\id_{I_{Z^\ast}(Y)}) \,.
  \ee
Evaluating this map explicitly on generic elements in $X$, $Y\,{\otimes^A}\, Z$ and 
$Z^\ast$, respectively, finally leads again to the expression given for
$\partial^\ihl_{X,Y,Z}$ in \eqref{eq:distformula}.
 \\[2pt]
The formulas for $\partial^\ihr_{X,Y,Z}$, $\widetilde{\partial}^\ihr_{X,Y,Z}$
and $\tensor*[^\ihr]{\delta}{^{Z^\ast}_X_Y}$ follow by analogous arguments.
\end{proof}

We also note the following properties of the distributors,
which follow independently of Lemma \ref{lem:iHom-strong}:

\begin{proposition}
Let $X$ be an \((A,A)\)-bimodule.
 \begin{enumerate} \addtolength{\itemsep}{1.8pt}
  \item \label{itm:linj}
If \(X\) is left \(\ota\)-flat, then
\(\partial^\ihl_{X,Y,Z}\) is injective for all \(Y,Z \iN A\bimod\).
  \item \label{itm:rsurj}
If \(X\) is right \(\oat\)-flat, then
\(\partial^\ihl_{Y,Z,X}\) is surjective for all \(Y,Z \iN A\bimod\).
   \item \label{itm:rinj}
If \(X\) is right \(\ota\)-flat, then
\(\partial^\ihr_{Y,Z,X}\) is injective for all \(Y,Z \iN A\bimod\).
   \item \label{itm:lsurj}
If \(X\) is left \(\oat\)-flat, then
\(\partial^\ihr_{X,Y,Z}\) is surjective for all \(Y,Z \iN A\bimod\).
  \end{enumerate}
\end{proposition}

\begin{proof}
Part \eqref{itm:linj}: Consider the second of 
the diagrams \eqref{eq:def:widetildepartial-l}. The composite
$\widetilde\gamma^{\ihl}_{X,Y,Z} \cir (X \otA \imath^{}_{Y,Z})$ 
is cobalanced, and \(\partial^\ihl_{X,Y,Z}\) is the unique homomorphism
determined by the universal property of the equalizer. Thus the kernel of
\(\partial^\ihl_{X,Y,Z}\) equals the kernel of $X\otA\imath^{}_{Y,Z}$.
Now note that this map is obtained by applying the functor \(X\otA{-}\) to the
inclusion in the equalizer definition of \(\oat\). Hence the kernel of
$X\otA\imath^{}_{Y,Z}$ is determined by torsion, and hence
\(\partial^l_{X,Y,Z}\) is injective if \(X\) is left \(\ota\)-flat.
 \\[2pt]
Part \eqref{itm:rsurj}: Recall that in the diagram \eqref{eq:square1b}
the composite $(\pi^{}_{X,Y} \oaT Z) \cir \gamma^{\ihl}_{X,Y,Z}$ 
is balanced, and that \(\partial^\ihl_{X,Y,Z}\) is the unique homomorphism determined
by the universal property of the coequalizer. Thus the image of $\pi^{}_{X,Y} \oaT Z$
equals the image of \(\partial^\ihl_{X,Y,Z}\). Note further that $\pi^{}_{X,Y} \oaT Z$
is obtained by applying the functor \(-\oat\) to the surjection in the coequalizer
definition of \(\ota\). Hence the image of \(\partial^\ihl_{X,Y,Z}\) is
determined by cotorsion, and thus \(\partial^\ihl_{X,Y,Z}\) is surjective if \(Z\)
is right \(\otimes^A\)-flat.
 \\[2pt]
The statements \eqref{itm:rinj} and \eqref{itm:lsurj} are shown analogously.
\end{proof}

It is worth pointing out that the failure of the distributors to be isomorphisms
can concern both their kernel and their image. We mention two concrete examples:
First, the three-dimensional algebra $\AA_3 \eq \ko[x,y]/{\langle x^2,y^2,xy\rangle}$
considered in \eqref{eq:def:A3} has two non-isomorphic three-dimensional
indecomposable modules, namely the projective $P \eq \AA_3 \eq \mathrm{span}_\ko\{1,x,y\}$
and the injective $I \eq P^\ast \eq \mathrm{span}_\ko\{1^\ast,x^\ast,y^\ast\}$.
We find that
  \be
  \begin{aligned}
  \mathrm{ker}(\partial^\ihl_{I,P,P}) & = \mathrm{span}_\ko
  \{\, [x^\ast\otik(x\otik x)] \,{-}\, [y^\ast\otik(y\otik x)] \,,
  \\
  & \phantom{= \mathrm{span}_\ko \{~ }
  [x^\ast\otik(x\otik y)] \,{-}\, [y^\ast\otik(y\otik y)] \,,
  [x^\ast\otik(y\otik x)] \,, 
  \\
  & \phantom{= \mathrm{span}_\ko \{~ }
  [x^\ast\otik(y\otik y)] \,, [y^\ast\otik(x\otik x)] \,, [y^\ast\otik(x\otik y)] \,\}\;,
  \nxl1
  \mathrm{im}(\partial^\ihl_{I,P,P}) & = \mathrm{span}_\ko
  \{\, [1^\ast\otik 1]\otik x,[1^\ast\otik 1]\otik y\,\} \,,
  \end{aligned}
  \ee
i.e.\ \(\partial^\ihl_{I,P,P}\) is neither injective nor surjective. 
Thus \(I\) is not \(\ota\)-flat, while \(P\) is not \(\oat\)-flat.
Second, for the algebra \eqref{eq:def:A2} of dual numbers, which up to isomorphism
has two indecomposable modules $S \eq \mathrm{span}_\ko\{s\}$ and
$P \eq A_2 \eq \mathrm{span}_\ko\{1,x\}$, we find an example for which the 
distributor vanishes: $\partial^\ihl_{S,P,S} \eq 0$.
So \(S\) is neither \(\ota\)-flat nor \(\oat\)-flat.

Both the concrete examples and the general results studied in this note 
reveal that GV-dualities are natural structures and that they deserve further 
study, in particular in the realm of quantum algebra and quantum topology.


\vfill

\noindent
{\sc Acknowledgments:}\\[.3em]
We thank Chelsea Walton for helpful comments.
J.F.\ is supported by VR under project no.\ 2022-02931.
C.S.\ is supported by the Deutsche Forschungsgemeinschaft (DFG, German Research 
Foundation) under SCHW1162/6-1 and under Germany's Excellence Strategy - EXC 2121 
``Quantum Universe'' - 390833306.
S.W.\ is supported by the Engineering and Physical Sciences Research Council (EPSRC)
EP/V053787/1 and by the Alexander von Humboldt Foundation.

\newpage


\begin{thebibliography}{ALSW}

\newcommand\wb{\,\linebreak[0]} \def\wB {$\,$\wb}
\newcommand\Bi[2]    {\bibitem[#2]{#1}}
\newcommand\COLL[2]  {{\em #1\/} ({#2})}
\newcommand\inBO[9]  {{\em #9}, in:\ {\em #1}, {#2}\ ({#3}, {#4} {#5}), p.\ {#6--#7} {\tt [#8]}}
\newcommand\J[7]     {{\em #7}, {#1} {#2} ({#3}) {#4--#5} {{\tt [#6]}}}
\newcommand\JO[6]    {{\em #6}, {#1} {#2} ({#3}) {#4--#5} }
\newcommand\BOOK[4]  {{\em #1\/} ({#2}, {#3} {#4})}
\newcommand\EBOO[2]  {{\sl #2}, available at {\tt http{/$\!$/#1}}}
\newcommand\Prep[2]  {{\em #2}, preprint {\tt #1}}

   \def\blms  {Bull.\wB London\wB Math.\wb Soc.}
   \def\coma  {Con\-temp.\wb Math.}
   \def\comp  {Com\-mun.\wb Math.\wb Phys.}
   \def\ijmb  {Int.\wb J.\wb Mod.\wb Phys.\ B}
   \def\izma  {Izvestiya: Math.}
   \def\jfst  {J.\wb Fac.\wb Sci.\wb Univ.\wB Tokyo}
   \def\joal  {J.\wB Al\-ge\-bra}
   \def\jpaa  {J.\wB Pure\wB Appl.\wb Alg.}
   \def\joms  {J.\wb Math.\wb Sci.}  
   \def\nupb  {Nucl.\wb Phys.\ B}
   \def\quto  {Quantum Topology}
   \def\taac  {Theo\-ry\wB and\wB Appl.\wb Cat.}
   \def\tams  {Trans.\wb Amer.\wb Math.\wb Soc.}
   \def\thcs  {Theor.\wb Computer\wB Science}

\Bi{alsW}  {ALSW} {R.\ Allen, S.\ Lentner, C.\ Schweigert, and S.\ Wood,
           \Prep{math.QA/2107.%
		 \linebreak[0]%
		05718} {Duality structures for module
           categories of vertex operator algebras and the Feigin Fuchs boson}}

\Bi{BArr}  {Ba1} {M.\ Barr, \BOOK{${}^*$-Autonomous Categories}
           {Springer Verlag}{Berlin}{1979}}
\Bi{barr5} {Ba2} {M.\ Barr, \JO\thcs{139}{1995}{115}{130}
           {Nonsymmetric *-autonomous categories}}
\Bi{barr8} {Ba3} {M.\ Barr, \JO\jpaa{111}{1996}{1}{20}
           {${}^*$-Autonomous categories, revisited}}
\Bi{barr9} {Ba4} {M.\ Barr, \JO\taac{6}{1999}{5}{24}
           {$\ast$-autonomous categories: once more around the track}}

\Bi{BOur3} {Bo} {N.\ Bourbaki, \BOOK{Algebra I $($Chapters 1-3$)$}
           {Springer Verlag}{Berlin}{1989}}

\Bi{boDr}  {BoD} {M.\ Boyarchenko and V.\ Drinfeld,
           \J\quto{4}{2013}{447}{489} {math.QA/1108.6020}
           {A duality formalism in the spirit of Grothendieck and Verdier}}

\Bi{brWo}  {BrW} {A.\ Brochier and L.\ Woike, \Prep{math.QA/2212.11259}
           {A classification of modular functors via factorization homology}}

\Bi{BRwi}  {BW} {T.\ Brzezi\'nski and R.\ Wisbauer,
           \BOOK{Corings and Comodules}
           {Cambridge University Press}{Cambridge}{2003}}

\Bi{coSe3} {CS} {J.R.B.\ Cockett and R.A.G.\ Seely, \JO\jpaa{114}{1997}{133}{173}
           {Weakly distributive categories}}

\Bi{EGno}  {EGNO} {P.I.\ Etingof, S.\ Gelaki, D.\ Nikshych, and V.\ Ostrik,
           \BOOK{Tensor Categories} {American Mathematical Society}{Providence}{2015}}

\Bi{FRbe2} {FB} {E.\ Frenkel and D.\ Ben-Zvi,
           \BOOK{Vertex Algebras and Algebraic Curves{\rm, second edition}}
           {American Mathematical Society}{Providence}{2004}}

\Bi{fScS2} {FSS} {J.\ Fuchs, G.\ Schaumann, and C.\ Schweigert, \J\tams{373}{2020}{1}{40}
           {math.RT/1612.04561}
           {Eilenberg-Watts calculus for finite categories and a bimodule Radford $S^4$ theorem}}

\Bi{gaKa}  {GK} {M.R.\ Gaberdiel and H.G.\ Kausch, \J\nupb{477}{1996}{293}{318}
           {hep-th/9604026} {Indecomposable fusion products}}

\Bi{garW3} {GRW} {M.R.\ Gaberdiel, I.\ Runkel, and S.\ Wood,
           \inBO{Conformal Field Theories and Tensor Categories}
           {C.M.\ Bai, J.\ Fuchs, Y.-Z.\ Huang, L.\ Kong, I.\ Runkel, and
           C.\ Schweigert, eds.} {Springer Verlag}{Berlin}{2014} {93}{168} {hep-th/1201.6273}
           {Logarithmic bulk and boundary conformal field theory and the full centre construction}}

\Bi{GAroz} {GR} {D.\ Gaitsgory and N.\ Rozenblyum,
           \BOOK{A Study in Derived Algebraic Geometry Volume I: Correspondences and Duality}
           {American Mathematical Society}{Providence}{2019}}

\Bi{haZo2} {HalZ} {S.\ Halbig and T.\ Zorman.  \Prep{math.CT/2301. 
                        \linebreak[0]
	   03545} {Duality in monoidal categories}}

\Bi{haLem} {HasL} {K.\ Hasegawa and L.\ Lemay, \JO\taac{33}{2018}{1145}{1157}
           {Linear distributivity with negation, star-autonomy, and Hopf monads}}

\Bi{HLZ}   {HLZ} {Y.-Z,\ Huang, J.\ Lepowsky, and L.\ Zhang,
           \inBO{Conformal Field Theories and Tensor Categories}
           {C.M.\ Bai et al., eds.} {Springer Verlag}{Berlin}{2014}
           {169}{248} {math.QA/1012.4193} {Logarithmic tensor product theory I}}
         
\Bi{ivanS} {Iv} {S.O.\ Ivanov, \JO\joms{183}{2012}{675}{680}
           {Nakayama functors and Eilenberg-Watts theorems}}

\Bi{kaRi}  {KR} {S.\ Kanade and D.\ Ridout,
           \inBO{Affine, Vertex and W-algebras}
           {D.\  Adamovi\'c and P.\ Papi, eds.}{Springer}{Cham}{2019} {135}{181}
           {math-ph/1812.10713}
           {NGK and HLZ: fusion for physicists and mathematicians}}

\Bi{lizu4} {LZ} {B.H.\ Lian and G.J.\ Zuckerman, \J\comp{154}{1993}{613}{646}
           {hep-th/9211072}
           {New perspectives on the BRST-algebraic structure of string theory}}

\Bi{mani27}{Ma} {Yu.I.\ Manin, \J\izma{81}{2017}{818}{826} {math.QA/1701.01261}
           {Grothendieck-Verdier duality patterns in quantum algebra}}

\Bi{muWo3} {MW} {L.\ M\"uller and L.\ Woike, \J\blms{53}{2021}{392}{403}
           {math.QA/2004.04689}
           {Dimensional reduction, extended topological field theories and orbifoldization}}

\Bi{nahm26}{Na} {W.\ Nahm, \J\ijmb{8}{1994}{3693}{3702}
           {hep-th/9402039} {Quasi-rational fusion products}}

\Bi{ncat}  {nLab} {nLab, {\em Dualizable module},
	   {\tt https:/$\!$/ncatlab.org/nlab/show/dualizable+module}}

\Bi{past}  {Pa} {C.\ Pastro, \J\taac{26}{2012}{194} {203} {math.CT/1010.5304}
           {Note on star-autonomous comonads}}

\Bi{seely3}{Se} {R.A.G.\ Seely, \JO\coma{92}{1989}{371}{382}
           {Linear logic, *-autonomous categories and cofree coalgebras}}

\Bi{shimi20}{Sh} {K.\ Shimizu, \J\joal{634}{2023}{237}{305} {math.CT/1904.00376}
           {Relative Serre functor for comodule algebras}}

\Bi{takeu3}{Ta} {M.\ Takeuchi, \JO\jfst{I\,A\,24}{1977}{629}{644}
           {Morita theorems for categories of comodules}}

\end{thebibliography}
\end{document}